\documentclass{article}
\usepackage{etex}

\usepackage{amsmath,amsthm,amssymb,amscd,mathtools}
\usepackage{amsfonts}
\usepackage{rotating}
\usepackage{euscript}
\usepackage{pst-node}
\usepackage{epsfig,verbatim}
\usepackage{pictexwd,dcpic}
\usepackage{enumerate}
\usepackage{caption}
\usepackage{verbatim}
\usepackage{float}

\def\be{\begin{equation}}
\def\ee{\end{equation}}

\def\C{{\mathbb C}} 
\def\f{\EuScript}
\def\N{{\mathbb N}} 
\def\P{{\mathbb P}}
\def\Z{{\mathbb Z}}

\def\ord{{\rm ord}}

\def\e{\eqref}
\def\phi{{\varphi}}
\def\v{{\varepsilon}} 
\def\tt{\widetilde}
\def\deg{{\rm deg\,}}

\def\cos{{\rm cos\,}}

\def\GCD{{\rm GCD }}
\def\LCM{{\rm LCM }}

\def\mod{{\rm mod\ }}

\def\bp{\begin{proposition}}
\def\ep{\end{proposition}}

\def\bt{\begin{theorem}}
\def\et{\end{theorem}}
\def\br{\begin{remark}}
\def\er{\end{remark}}
\def\be{\begin{equation}}
\def\bee{\begin{equation*}}
\def\l{\label}

\def\e{\eqref}

\def\ee{\end{equation}}
\def\eee{\end{equation*}}
\def\bl{\begin{lemma}}
\def\el{\end{lemma}}
\def\bc{\begin{corollary}}
\def\ec{\end{corollary}}
\def\pr{\noindent{\it Proof. }}

\def\bd{\begin{definition}}
\def\ed{\end{definition}}
\def\t{\widetilde}

\newtheorem{theorem}{Theorem}[section]
\newtheorem{lemma}[theorem]{Lemma}
\newtheorem{corollary}[theorem]{Corollary}
\newtheorem{proposition}[theorem]{Proposition}
\newtheorem{definition}[theorem]{Definition}
\newtheorem{remark}[theorem]{Remark}

\mathtoolsset{showonlyrefs}

\begin{document}
\title{On generalized Latt\`es maps}
\author{F. Pakovich}
\maketitle
\begin{abstract}
We introduce a class of rational functions $A:\,\C\P^1\rightarrow \C\P^1$
which can be considered as a natural 
extension of the class of Latt\`es maps, and establish basic properties of functions from this class.

\end{abstract}

\begin{section}{Introduction}

Latt\`es maps are rational functions $A:\,\C\P^1\rightarrow \C\P^1$ of degree at least two which can be characterized in one of the following equivalent ways (see \cite{mil2}). 
First, a Latt\`es map $A$ can be defined by the condition that
there exist a compact Riemann surface $R$ of genus one  
 and holomorphic maps $B: R\rightarrow R$ and \linebreak $\pi:R\rightarrow \C\P^1$
such that the diagram 
\be \l{xxuuii}
\begin{CD}
R @>B>> R \\
@VV\pi V @VV\pi V\\ 
\C\P^1 @>A >> \ \ \C \P^1
\end{CD}
\ee
commutes. 
This condition can be replaced by
the apparently stronger condi\-tion that there exists a  diagram as above such that  $\pi$  is the quotient map
$\pi: {R}\rightarrow   
{R}/ \Gamma$ for some finite subgroup  $ \Gamma$ of the automorphism group ${\rm Aut}({R})$. 
Finally, Latt\`es maps  can be  characterized  in terms of their ramification.

The last characterization uses the notion of orbifold. 
By definition, an {\it orbifold} $\f O$ on $\C\P^1$  is a ramification function $\nu:\C\P^1\rightarrow \mathbb N$ which takes the value $\nu(z)=1$ except at a finite set of points. We always will assume that considered orbifolds are {\it good} meaning that we forbid $\f O$ to have exactly one point with $\nu(z)\neq 1$ or  two such points $z_1,$ $z_2$ with $\nu(z_1)\neq \nu(z_2).$ 
A rational function $f$
is called  {\it a covering map} $f:\,  \f O_1\rightarrow \f O_2$
between orbifolds with ramifications functions $\nu_1$ and $\nu_2$ 
if for any $z\in \C\P^1$ the equality 
\be \l{u} \nu_{2}(f(z))=\nu_{1}(z)\deg_zf\ee holds.
In these terms, a Latt\`es map can be defined as a rational function $A$  such that $A:\f O\rightarrow \f O$ is a  covering {\it self-map}
for some orbifold $\f O$.

In the recent paper \cite{semi} a class of rational functions $A$  satisfying  \eqref{xxuuii} under the assumption that the surface
$R$ is the {\it Riemann sphere} was  considered. 
It was shown in \cite{semi}  that under certain  restrictions such functions
 posses a number of remarkable properties 
similar to properties of Latt\`es maps. In particular, they are related to finite subgroups of the group $Aut(\C\P^1)$, and admit a description 
in terms of orbifolds. In this paper, modifying the approach of \cite{semi}, we construct 
 a unified theory
which equally fits the classical Latt\`es maps and functions studied in \cite{semi}, using the term ``generalized  Latt\`es maps" for the set of functions obtained in this way.

Notice that allowing $R$ in \eqref{xxuuii} to be an {\it arbitrary} compact Riemann surface does not lead to a yet more general class of functions,  
since for $R$ of genus at least two any holomorphic map $B: R\rightarrow R$ has degree one. 
Notice also that in order to define an interesting class of functions $A$ through diagram \eqref{xxuuii} 
with $R=\C\P^1$ some restrictions on $A$, $B$, and $\pi$ are necessary, 
since  there exist 
too many  rational functions making diagram \eqref{xxuuii} commutative. Say, for any rational functions $U$ and $V$ the diagram 
\be 
\begin{CD}
\C\P^1 @>U\circ V>> \C\P^1 \\
@VV {V} V @VV {V} V\\ 
\C\P^1 @>V\circ U >> \ \ \C \P^1
\end{CD}
\ee
commutes, and it is clear that the function $V\circ U$ does not posses any special properties in general.

 The easiest way to define generalized Latt\`es maps uses the concept of a  minimal holomorphic  map
between orbifolds. By definition, a rational function $f$ is  called {\it a  minimal holomorphic  map}  $f:\,  \f O_1\rightarrow \f O_2$ 
between orbifolds
if  for any $z\in \C\P^1$
the condition 
$$ \nu_{2}(f(z))=\nu_{1}(z)\GCD(\deg_zf, \nu_{2}(f(z))$$ holds.
It is easy to see that any  covering map $A:\f O_1\rightarrow \f O_2$  between orbifolds is a minimal holomorphic  map, but  the inverse is not true. Say that a rational function $A$ of degree at least two is a {\it generalized Latt\`es map} if there exists an orbifold $\f O$ distinct from the non-ramified sphere such that  $A:\,  \f O\rightarrow \f O$ is a  minimal holomorphic map between orbifolds. 
 
Recall that for an orbifold $\f O$ the {\it  Euler characteristic} of $\f O$ is the number
$$ \chi(\f O)=2+\sum_{z\in \C\P^1}\left(\frac{1}{\nu(z)}-1\right),$$
the set of {\it singular points} of $\f O$ is the set 
$$c(\f O)=\{z_1,z_2, \dots, z_s, \dots \}=\{z\in \C\P^1 \mid \nu(z)>1\},$$ and  the {\it signature} of $\f O$ is the set 
$$\nu(\f O)=\{\nu(z_1),\nu(z_2), \dots , \nu(z_s), \dots \}.$$ 
It is well-known that if  $A:\f O\rightarrow \f O$ is a covering map between orbifolds, then 
  the Euler characteristic  of  $\f O$
 equals zero, implying  that
 the signature of $\f O$ belongs to the list \be \l{list} \{2,2,2,2\}, \ \ \ \{3,3,3\}, \ \ \ \{2,4,4\}, \ \ \ \{2,3,6\}.\ee
On the other hand, if  $A:\f O\rightarrow \f O$ is a  minimal holomorphic  map  between orbifolds, then the Euler characteristic of $\f O$ is {\it non-negative}. Thus,  to the above list we should add the signatures \be \l{list2} \{n,n\}, \ \ n\geq 2,\ \ \ \ \{2,2,n\}, 
\ \ n\geq 2,\ \ \ \ \{2,3,3\}, \ \ \ \  \{2,3,4\}, \ \ \  \ \{2,3,5\},\ee corresponding to orbifolds of positive Euler characteristic.

In this paper we provide three characterizations of   generalized Latt\`es maps parallel to three  characterizations of  Latt\`es maps given in the paper \cite{mil2}  by J. Milnor. 
Let $R_1,$ $R_2$, and $R'$ be Riemann surfaces. Say that a holomorphic map  $h:\, R_1\rightarrow R'$ 
 is 
a {\it compositional right factor} of a holomorphic map  
$f:\, R_1\rightarrow R_2$ if there exists a holomorphic map $g:\, R'\rightarrow R_2$
 such that $f=g\circ h$. Compositional left factors are defined similarly.  In this notation, the
following statement holds.

\bt \l{bas} Let $A$ be a rational function of degree at least two. Then the following conditions are equivalent.
\begin{enumerate}

\item There exist a compact Riemann surface $R$  of genus zero or one 
and holomorphic maps $B: R\rightarrow R$  and $\pi:R\rightarrow \C\P^1$
such that  
the diagram 
\be \l{xui}
\begin{CD}
R @>B>> R \\
@VV\pi V @VV\pi V\\ 
\C\P^1 @>A >> \ \ \C \P^1
\end{CD}
\ee
commutes, and $\pi$ is not a compositional right  factor of $B^{\circ s}$ for some $s\geq 1.$  

\item There exist a compact Riemann surface $ {R}$ of genus zero or one, a finite non-trivial group $\Gamma\subseteq Aut( {R})$, an isomorphism $\phi:\, \Gamma \rightarrow \Gamma $, 
and  a holomorphic map $B:\, {R} \rightarrow   {R}$   
such that the diagram 
\be \l{xui+}
\begin{CD}
 {R} @>  B>>  {R} \\
@VV  \pi V @VV \pi  V\\ 
\C\P^1 @>A >> \ \ \C \P^1\, ,
\end{CD}
\ee
where $\pi:\, R\rightarrow R/\Gamma$ is the quotient map,  commutes, and for any $\sigma\in \Gamma$ the equality
\be \l{homm2}  B\circ\sigma=\phi(\sigma)\circ B \ee holds. 
\item There exists an orbifold $\f O$, distinct from the non-ramified sphere, such that  $$A:\,  \f O\rightarrow \f O$$ is a  minimal holomorphic map between orbifolds. 

\end{enumerate}

\et

Let us make several comments concerning conditions of Theorem \ref{bas}. By definition, 
$A:\,  \f O\rightarrow \f O$ is a  minimal holomorphic map between orbifolds if 
\be \l{uu} \nu(A(z))=\nu(z)\GCD(\deg_zA,\nu (A(z))), \ \ \ \ \ \ \ \\ z\in \C\P^1,\ee and 
it is easy to see that for the Riemann sphere, considered as a non-ramified orbifold, 
this condition holds for {\it any} rational function $A$. Thus, we must exclude this case in the third condition. 
By the same reason, we assume that $\Gamma\neq \{e\}$ 
in the second condition. 

The assumption in the first condition, requiring that $\pi$ is not a compositional right  factor of some iterate of $B$,  
is always satisfied if $g(R)=1$,   since for any decomposition
$$R\overset{\pi}{\longrightarrow}R^{\prime}\overset{w}{\longrightarrow}R$$
of $B^{\circ s}$, $s\geq 1,$  the genus of $R^{\prime}$ must be equal to one.
However, this assumption is essential if $R=\C\P^1.$
It can be replaced by the assumption  that $\pi$ is not a compositional {\it left}  factor of some iterate of $A$.  
Further, notice that for any diagram \eqref{xui+} such that  $\pi:\, R\rightarrow R/\Gamma$ is the quotient map for some finite group $\Gamma\subseteq Aut( {R})$, condition \eqref{homm2}
holds for some {\it homomorphism}  $\phi:\, \Gamma \rightarrow \Gamma $. Moreover, this homomorphism is always an isomorphism if $g(R)=1$,
however may have a non-trivial kernel if $R=\C\P^1$.

The paper is organized as follows. 
In the second section we recall main technical results of \cite{semi} about Riemann surfaces orbifolds and different kinds of maps between orbifolds. 
In the  third section we describe 
a general structure of holomorphic maps satisfying the semiconjugacy condition \eqref{xxuuii}, where $R$ is 
a compact Riemann surface  of genus zero or one, and prove Theorem \ref{bas}. 
In the fourth section we study properties of  generalized Latt\`es maps related to the operations of composition 
and decomposition. 
In the 
fifth section 
 we describe rational functions satisfying condition \eqref{uu} for orbifolds $\f O$ with signatures $\{n,n\}$, $n\geq 2,$ and  $\{2,2,n\},$ $n> 2.$ 

In the sixth  section we investigate the following problem: given a rational function $A$, 
what are orbifolds
$\f O$ such that   $A:\,  \f O\rightarrow \f O$ is a  minimal holomorphic  map between orbifolds ? For ordinary Latt\`es maps, there exists at most one such an orbifold defined by dynamical properties of $A$. On the other hand,
for generalized Latt\`es maps there might be several and even infinitely many such orbifolds.  For example, it is easy to see that $z^{\pm n}:\f O\rightarrow \f O$ is a minimal holomorphic map for any $\f O$ defined  by  $$\nu(0)=m, \ \ \ \ \nu(\infty)=m, \ \ \ \GCD(n,m)=1,$$ while
 $\pm T_{n}:\f O\rightarrow \f O$ is a minimal holomorphic map  for any $\f O$ defined by  the conditions $$
\nu(-1)=\nu(1)=2, \ \ \ \nu(\infty)=m, \ \ \ \  \GCD(n,m)=1.$$
Nevertheless, we show  that if $A$ is  not conjugate to $z^{\pm n}$ or $\pm T_n$, then there exists  
a ``maximal'' orbifold $\f O$ such that \eqref{uu} holds. 
In more details, 
for   orbifolds $\f O_1$  and $\f O_2$ 
write $ \f O_1\preceq \f O_2$ if for any $z\in\C\P^1$
the condition   $\nu_1(z)\mid \nu_2(z)$ holds.
In this notation, the main result of the sixth section and one of the main results of the paper is following.

\bt \l{uni} 
Let $A$ be a rational function  of degree at least two not conjugate to $z^{\pm d}$ or $\pm T_d.$ Then there exists an orbifold $\f O_0^A$ such that $A:\, \f O_0^A\rightarrow \f O_0^A$
is a minimal holomorphic map between orbifolds, and for any orbifold $\f O$ such that 
$A:\, \f O\rightarrow \f O$ is a minimal holomorphic map between orbifolds the relation $\f O\preceq \f O_0^A$ holds. Furthermore, $\f O_0^{A^{\circ l}}=\f O_0^A$ for any 
$l\geq 1$.
\et

In the seventh section we relate the problem of describing 
generalized Latt\`es maps
which are not ordinary Latt\`es maps with the problem of
describing rational functions commuting with a finite automorphism group  of $Aut(\C\P^1).$ 
We recall a description of such functions obtained by Doyle and McMullen (\cite{dm}), and give examples of practical calculations of corresponding generalized Latt\`es maps of small degrees.  
Finally, we show that 
{\it polynomial} generalized Latt\`es maps reduce to the  series $T_n$ and $z^rR^n(z),$ where $R\in \C[z]$ and $\GCD(r,n)=1,$ emerging in the Ritt theory of polynomial decompositions \cite{r2}. 

\end{section}

\begin{section}{Orbifolds and maps between orbifolds}
In this section we recall basic definitions concerning Riemann surface orbifolds 
 (see \cite{mil}, Appendix E), 
and overview 
some   technical results  obtained in the paper \cite{semi}.

A Riemann surface orbifold is a pair $\f O=(R,\nu)$ consisting of a Riemann surface $R$ and a ramification function $\nu:R\rightarrow \mathbb N$ which takes the value $\nu(z)=1$ except at isolated points. 
For an orbifold $\f O=(R,\nu)$ 
 the {\it  Euler characteristic} of $\f O$ is the number
$$ \chi(\f O)=\chi(R)+\sum_{z\in R}\left(\frac{1}{\nu(z)}-1\right),$$
the set of {\it singular points} of $\f O$ is the set 
$$c(\f O)=\{z_1,z_2, \dots, z_s, \dots \}=\{z\in R \mid \nu(z)>1\},$$ and  the {\it signature} of $\f O$ is the set 
$$\nu(\f O)=\{\nu(z_1),\nu(z_2), \dots , \nu(z_s), \dots \}.$$ 
For   orbifolds $\f O_1=(R_1,\nu_1)$  and $\f O_2=(R_2,\nu_2)$ 
write \be \l{elki} \f O_1\preceq \f O_2 \ee 
if $R_1=R_2$, and for any $z\in R_1$ the condition $$\nu_1(z)\mid \nu_2(z)$$ holds.
Clearly, \eqref{elki} implies that 
$$\chi(\f O_1)\geq \chi(\f O_2).$$

Let $R_1$, $R_2$ be Riemann surfaces provided with ramification functions $\nu_1,$ $\nu_2$.
A holomorphic branched covering map
$f:\, R_1\rightarrow R_2$ 
is called  {\it a covering map} $f:\,  \f O_1\rightarrow \f O_2$
between orbifolds
$\f O_1=(R_1,\nu_1)$ and $\f O_2=(R_2,\nu_2)$
if for any $z\in R_1$ the equality 
\be \l{us} \nu_{2}(f(z))=\nu_{1}(z)\deg_zf\ee holds, where $\deg_zf$ is the local degree of $f$ at the point $z$.
If for any $z\in R_1$ instead of equality \eqref{us} 
a weaker condition 
\be \l{uuss} \nu_{2}(f(z))\mid \nu_{1}(z)\deg_zf\ee
holds,  then $f$
is called {\it a holomorphic map} $f:\,  \f O_1\rightarrow \f O_2$
between orbifolds
$\f O_1$ and $\f O_2.$

{\it A universal covering} of an orbifold ${\f O}$
is a covering map between orbifolds \linebreak $\theta_{\f O}:\,
\tt {\f O}\rightarrow \f O$ such that $\tt R$ is simply connected and $\tt \nu(z)\equiv 1.$ 
If $\theta_{\f O}$ is such a map, then 
there exists a group $\Gamma_{\f O}$ of conformal automorphisms of $\tt R$ such that the equality 
$\theta_{\f O}(z_1)=\theta_{\f O}(z_2)$ holds for $z_1,z_2\in \tt R$ if and only if $z_1=\sigma(z_2)$ for some $\sigma\in \Gamma_{\f O}.$ A universal covering exists and 
is unique up to a conformal isomorphism of $\tt R,$
unless $\f O$ is the Riemann sphere with one ramified point or with two ramified points $z_1,$ $z_2$ such that $\nu(z_1)\neq \nu(z_2).$   
 Furthermore, 
$\tt R=\mathbb D$ if and only if $\chi(\f O)<0,$ $\tt R=\C$ if and only if $\chi(\f O)=0,$ and $\tt R=\C\P^1$ if and only if $\chi(\f O)>0$ (see e. g. \cite{fk}, Section IV.9.12).
Abusing  notation we will use the symbol $\tt {\f O}$ both for the
orbifold and for the  Riemann surface  $\tt R$.

Covering maps between orbifolds lift to isomorphisms between their universal co\-ve\-rings.
More generally, for holomorphic maps the following proposition holds (see \cite{semi}, Propo\-sition 3.1).

\bp \l{poiu} Let $f:\,  \f O_1\rightarrow \f O_2$ be a holomorphic map between orbifolds. Then for any choice of $\theta_{\f O_1}$ and $\theta_{\f O_2}$ there exist 
a holomorphic map $F:\, \tt {\f O_1} \rightarrow \tt {\f O_2}$ and 
a homomorphism $\phi:\, \Gamma_{\f O_1}\rightarrow \Gamma_{\f O_2}$ such that the diagram 
\be \l{dia2}
\begin{CD}
\tt {\f O_1} @>F>> \tt {\f O_2}\\
@VV\theta_{\f O_1}V @VV\theta_{\f O_2}V\\ 
\f O_1 @>f >> \f O_2\ 
\end{CD}
\ee
is commutative and 
for any $\sigma\in \Gamma_{\f O_1}$ the equality
\be \l{homm}  F\circ\sigma=\phi(\sigma)\circ F \ee holds.
The map $F$ is defined by $\theta_{\f O_1}$, $\theta_{\f O_2}$, and $f$  
uniquely up to a transformation 
$F\rightarrow g\circ F,$ where $g\in \Gamma_{\f O_2}$. 
In the other direction, for any holomorphic map \linebreak $F:\, \tt {\f O_1} \rightarrow \tt {\f O_2}$  which satisfies \eqref{homm} for some homomorphism $\phi:\, \Gamma_{\f O_1}\rightarrow \Gamma_{\f O_2}$
there exists a uniquely defined  holomorphic map between orbifolds $f:\,  \f O_1\rightarrow \f O_2$ such that diagram \eqref{dia2} is commutative.
The holomorphic map $F$ is an isomorphism if and only if $f$ is a covering map between orbifolds. \qed

\ep

If $f:\,  \f O_1\rightarrow \f O_2$ is a covering map between orbifolds with compact $R_1$ and $R_2$, then  the Riemann-Hurwitz 
formula implies that 
\be \l{rhor} \chi(\f O_1)=d \chi(\f O_2), \ee
where $d=\deg f$. 
For holomorphic maps the following statement is true (see \cite{semi}, Proposition 3.2). 

\bp \l{p1} Let $f:\, \f O_1\rightarrow \f O_2$ be a holomorphic map between orbifolds with compact $R_1$ and $R_2$.
Then 
\be \l{iioopp} \chi(\f O_1)\leq \chi(\f O_2)\,\deg f, \ee and the equality 
holds if and only if $f:\, \f O_1\rightarrow \f O_2$ is a covering map between orbifolds. \qed
\ep

Let $R_1$, $R_2$ be Riemann surfaces and 
$f:\, R_1\rightarrow R_2$ a holomorphic branched covering map. Assume that $R_2$ is provided with ramification function $\nu_2$. In order to define a ramification function $\nu_1$ on $R_1$ so that $f$ would be a holomorphic map between orbifolds $\f O_1=(R_1,\nu_1)$ and $\f O_2=(R_2,\nu_2)$ 
we must satisfy condition \eqref{uuss}, and it is easy to see that
for any  $z\in R_1$ a minimal possible value for $\nu_1(z)$ is defined by 
the equality 
\be \l{rys} \nu_{2}(f(z))=\nu_{1}(z)\GCD(\deg_zf, \nu_{2}(f(z)).\ee 
In case if \eqref{rys} is satisfied for  any $z\in R_1$ we 
say that $f$ is {\it a  minimal holomorphic  map} 
between orbifolds 
$\f O_1=(R_1,\nu_1)$ and $\f O_2=(R_2,\nu_2)$.

It follows from the definition that for any orbifold $\f O=(R,\nu)$ and holomorphic branched covering map $f:\, R^{\prime} \rightarrow R$ there exists a unique orbifold structure $\nu^{\prime}$ on $R^{\prime}$ such that 
$f$ becomes a minimal holomorphic map between orbifolds. We will denote the corresponding orbifold by $f^*\f O.$
Notice that any covering map between orbifolds $f:\,  \f O_1\rightarrow \f O_2$ is a  minimal holomorphic map. In particular, $\f O_1=f^*\f O_2.$
For   orbifolds $\f O_1$  and $\f O_2$ we will 
write
\be \l{edc} \nu(\f O_1)\leq \nu(\f O_2)\ee
if 
for any $x\in c(\f O_1)$ there exists $y\in c(\f O_2)$ such that $\nu(x)\mid \nu(y).$ Clearly, the condition that $f:\,  \f O_1\rightarrow \f O_2$ is a  minimal holomorphic map 
 implies condition \eqref{edc}.  
Notice that \eqref{elki} implies \eqref{edc}
 but the inverse is not true in general.

\vskip 0.2cm
Minimal holomorphic maps between orbifolds possess the following fundamental property (see \cite{semi}, Theorem 4.1).

\bt \l{serrr} Let $f:\, R^{\prime\prime} \rightarrow R^{\prime}$ and $g:\, R^{\prime} \rightarrow R$ be holomorphic branched covering maps, and  $\f O=(R,\nu)$ an orbifold. 
Then 

$$(g\circ f)^*\f O= f^*(g^*\f O).\eqno{\Box}$$
\et

Theorem \ref{serrr} implies in particular the following corollaries (see   \cite{semi}, Corollary 4.1 and Corollary 4.2).

\bc \l{serka0} Let $f:\, \f O_1\rightarrow \f O^{\prime}$ and $g:\, \f O^{\prime}\rightarrow \f O_2$ be minimal holomorphic maps (resp. covering maps) between orbifolds.
Then  $g\circ f:\, \f O_1\rightarrow \f O_2$ is  a minimal holomorphic map (resp. covering map). \qed
\ec

\bc \l{indu2}  Let $f:\, R_1 \rightarrow R^{\prime}$ and $g:\, R^{\prime} \rightarrow R_2$ be holomorphic branched covering maps, and  $\f O_1=(R_1,\nu_1)$ and  $\f O_2=(R_2,\nu_2)$
orbifolds. Assume that  $g\circ f:\, \f O_1\rightarrow \f O_2$ is  a minimal holomorphic map (resp. a co\-vering map). Then $g:\, g^*\f O_2\rightarrow \f O_2$  and  $f:\, \f O_1\rightarrow g^*\f O_2 $ are minimal holomorphic maps (resp. covering maps). \qed
\ec 


With each holomorphic map $f:\, R_1\rightarrow R_2$ between compact Riemann surfaces 
one can associate in a natural way two orbifolds $\f O_1^f=(R_1,\nu_1^f)$ and 
$\f O_2^f=(R_2,\nu_2^f)$, setting $\nu_2^f(z)$  
equal to the least common multiple of local degrees of $f$ at the points 
of the preimage $f^{-1}\{z\}$, and $$\nu_1^f(z)=\nu_2^f(f(z))/\deg_zf.$$
By construction,  $f:\, \f O_1^f\rightarrow \f O_2^f$ 
is a covering map between orbifolds. 
It is easy to  see that the covering map $f:\, \f O_1^f\rightarrow \f O_2^f$ is minimal in the following sense. For any covering map between orbifolds $f:\, \f O_1\rightarrow \f O_2$ we have:
\be \l{elki+} \f O_1^f\preceq \f O_1, \ \ \ \f O_2^f\preceq \f O_2.\ee
On the other hand, for any holomorphic map $f:\, \f O_1\rightarrow \f O_2$ we 
have:
$$ f^{*} \f O_2 \preceq\f O_1.$$

\vskip 0.2cm

Orbifolds  $\f O_1^f$ and 
$\f O_2^f$ are useful for the study of the functional equation
\be \l{m} f\circ p=g\circ q, \ee 
where  
$$p:\, R\rightarrow C_1,  \ \ \ \ f:\, C_1\rightarrow \C\P^1,\ \ \ \ q:\, R\rightarrow C_2, \ \ \ \ g:\, C_2\rightarrow \C\P^1$$ 
are holomorphic maps between compact Riemann surfaces.
Recall that the fiber product  of the coverings $f:C_1\rightarrow \C\P^1$ and $g:C_2\rightarrow \C\P^1$
is defined as the set of pairs $(z_1,z_2)\in C_1\times C_2$ such that $f(z_1)=g(z_2).$ 
 The fiber product is a finite union of singular Riemann surfaces, and can be described in terms
of the monodromy groups of $f$ and $g$  (see e.g. \cite{pak},  Section 2).
Say that a solution $f,p,g,q$ of \eqref{m} is {\it good} if the fiber product of $f$ and $g$  consists of a unique component, 
and $p$ and $q$ have no {\it non-trivial common compositional right factor}. By definition, the last condition means that if
\be \l{kaban} p= \tt p\circ  w, \ \ \ q= \tt q\circ w\ee for some holomorphic maps $$w:\, R \rightarrow \tt R,  \ \ \ \ \tt p:\, \tt R\rightarrow C_1, \ \ \ \ \tt q:\, \tt R\rightarrow C_2,$$ then necessarily $\deg w=1.$
Notice that if $f$ and $g$ are  rational functions, then the fiber product of $f$ and $g$ has a unique component if and only if  the algebraic curve $$f(x)-g(y)=0$$
is irreducible. 
On the other hand, the L\"uroth theorem implies that if $p$ and $q$ are rational functions, then they  have no non-trivial common compositional right factor if and only if
$\C(p,q)=\C(z)$.

\vskip 0.2cm

In the above notation the following statement holds (see \cite{semi}, Theorem 4.2).

\bt \l{t1} Let $f,p,g,q$ be a good solution of \eqref{m}.
Then the commutative diagram 
\be 
\begin{CD}
\f O_1^q @>p>> \f O_1^f\\
@VV q V @VV f V\\ 
\f O_2^q @>g >> \f O_2^f\ 
\end{CD}
\ee
consists of minimal holomorphic  maps between orbifolds. \qed  
\et

Below we will use the following criterion (see \cite{semi}, Lemma 2.1).

 \bl \l{good} A 
solution $f,p,g,q$ of \e{m} is good whenever 
any two of the following three conditions are satisfied:

\begin{itemize}
\item the fiber product of $f$ and $g$ has a unique component,
\item $p$ and $q$ have no non-trivial common compositional right factor,
\item $ \deg f=\deg q, \ \ \ \deg g=\deg p.$  \qed
\end{itemize}
\el

In this paper essentially all considered orbifolds will be defined on $\C\P^1.$ 
The only exceptions from this rule are orbifolds which are universal coverings.
So, usually we will omit  the Riemann surface $R$ in the definition of $\f O=(R,\nu)$
meaning that $R=\C\P^1.$ We also will assume that all considered orbifolds have a universal covering.  

The central role in our exposition is played by orbifolds $\f O$ of non-negative Euler characteristic. 
For such orbifolds the corresponding groups $\Gamma_{\f O}$ and functions $\theta_{\f O}$ are described as follows.
Groups $\Gamma_{\f O}\subset Aut(\C)$ corresponding to orbifolds $\f O$ with signatures \eqref{list}  
are generated by translations of $\C$ by elements of some lattice $L\subset \C$ of rank two and the 
rotation $z\rightarrow  \v z,$ where $\v$ is an $n$th root of unity with $n$ equal to 2,3,4, or 6, such that  $\v L=L$.
In more details, the subgroup  $\Lambda_{\f O}\subset \Gamma_{\f O}$  generated by all translations
is a free group of rank two so that $R=\C/\Lambda_{\f O}$ is a torus, $\Lambda_{\f O}$ is normal in $\Gamma_{\f O}$, 
and   $\Gamma_{\f O}/\Lambda_{\f O}$ is a cyclic group of order 2,3,4, or 6, which acts as a group of automorphisms of  $R=\C/\Lambda_{\f O}$.
 Accordingly, the functions $\theta_{\f O}$ 
may be written in terms of the  corresponding
Weierstrass functions as $\wp(z),$ $\wp^{\prime }(z),$ $\wp^2(z),$  and $\wp^{\prime 2}(z)$  (see \cite{fk}, 
Section IV.9.5 and \cite{mil2}).  

Groups $\Gamma_{\f O}\subset Aut(\C\P^1)$ corresponding to   orbifolds $\f O$ with signatures \eqref{list2} are the well-known finite subgroups 
 $C_n,$  $D_{2n},$  $A_4,$ $S_4,$ $A_5$ of $Aut(\C\P^1)$, and the functions $\theta_{\f O}$ are Galois coverings of $\C\P^1$ by $\C\P^1$ of degrees 
$n$, $2n,$ $12,$ $24,$ $60,$ calculated for the first time by Klein in \cite{klein}.

\vskip 0.2cm
In conclusion of this section, let us mention the following  
more precise version of Proposition \ref{poiu} for minimal holomorphic self-maps between orbifolds of positive characteristic (see \cite{semi}, Theorem 5.1).

\bt \l{las}  Let $A$ and $F$ be rational functions of degree at least two and  $\f O$ an orbifold with $\chi(\f O)>0$ such that 
$A:\, \f O \rightarrow \f O$ 
is a holomorphic map between orbifolds  and  the diagram 
\be \l{dia3}
\begin{CD}
\tt {\f O} @>F>> \tt {\f O}\\
@VV\theta_{\f O}V @VV\theta_{\f O}V\\ 
\f O @>A >> \f O\ 
\end{CD}
\ee
commutes.
Then the following conditions are equivalent.

\begin{enumerate}
\item The holomorphic map $A$ is a minimal holomorphic  map. 
\item  The homomorphism $\phi:\, \Gamma_{\f O}\rightarrow \Gamma_{\f O}$ defined by the equality 
\be \l{homo} F\circ\sigma=\phi(\sigma)\circ F, \ \ \ \sigma\in \Gamma_{\f O},\ee is an automorphism of $\Gamma_{\f O}$.
\item The triple $F,$ $A,$ $\theta_{\f O}$ is a good solution of the equation 
$$ 
A\circ \theta_{\f O}=\theta_{\f O}\circ F. \eqno{\Box}$$ 
\end{enumerate}

\et

\end{section} 

\begin{section}{Semiconjugacies and generalized Latt\`es maps}

In this section we describe a general structure of holomorphic maps satisfying the semiconjugacy condition \eqref{xxuuii}, where $R$ is 
a compact Riemann surface  of genus zero or one, and prove Theorem \ref{bas}.  
Recall that we defined a {\it generalized Latt\`es map} as a rational function of degree at least two such that $A:\,  \f O\rightarrow \f O$ is a  minimal holomorphic map between orbifolds for some $\f O$ distinct from the non-ramified sphere. 
By Proposition \ref{p1}, for such $\f O$ necessarily $\chi(\f O)\geq 0.$ Notice that if $\chi(\f O)=0$, then $A:\,  \f O\rightarrow 
\f  O$ is a covering map
by Proposition \ref{p1}, and therefore $A$ is an ordinary Latt\`es map.

Let $B$ be a rational function of degree at least two. 
For any decomposition $B=V\circ U,$ where $U$ and $V$ are rational functions, the 
rational function $\t B=U\circ V$ is called an elementary transformation of $B$, and rational functions $B$ and $A$ are called  {\it equivalent}  if there exists 
a chain of elementary transformations between $B$ and $A$. For a rational function $B$ we will denote its equivalence class by $[B].$
Since for any invertible rational function $W$ the equality
$$B=(B\circ W)\circ W^{-1}$$ holds, each equivalence class $[B]$ is a union of conjugacy classes. 
Thus, the relation $\sim$ can be considered as a 
weaker form of the classical conjugacy relation.  Notice that 
an equivalence class $[B]$ 
contains infinitely many conjugacy classes if and only if 
$B$ is a flexible Latt\`es map (see \cite{rec}).

The  connection between the relation $\sim$ and semiconjugacy is straightforward. Namely, for $\t B$ and $B$ as above we have:  
$$\t B\circ U=U\circ B, \ \ \ \ \  B\circ V=V\circ \t B,$$ implying inductively that if
$B\sim \t B$, then $B$ is semiconjugate to $\t B$, and   $\t B$ is semiconjugate to $B.$ Moreover, the following statement is true.  

\bl \l{lem1} Let \be \l{chh} B\rightarrow B_1 \rightarrow B_2  \rightarrow \dots \rightarrow B_s\ee 
be a chain of elementary transformations, and   $U_i,$ $V_i,$ $1\leq i \leq s,$  rational functions such that 
$$B=V_1\circ U_1, \ \ \  B_i= U_i\circ V_i, \ \ \ \ \ 1\leq i\leq s,$$ 
and
\be \l{seqs}  U_{i}\circ V_{i}=V_{i+1}\circ U_{i+1},\ \ \ 1 \leq i \leq s-1.\ee
Then  the functions
$$U=U_s\circ U_{s-1}\circ \dots \circ U_{1}, \ \ \ \ V=V_{1}\circ \dots \circ V_{s-1}\circ V_s$$
make the diagram 
\be 
\begin{CD} 
\C\P^1 @> B>>\C\P^1 \\ 
@V U  VV @VV U   V\\ 
 \C\P^1 @> B_s>> \C\P^1
 \\ 
@V {V}  VV @VV {V}   V\\ 
 \C\P^1 @> B>> \C\P^1, 
\end{CD} 
\ee
commutative 
and 
satisfy the equalities 
$$V\circ U=B^{\circ s}, \ \ \ \ \  \ U\circ V=B_s^{\circ s}.$$ 
\el
\pr Indeed, we have:
$$
\begin{gathered}
B_s\circ (U_s\circ U_{s-1}\circ \dots \circ U_{1})=U_s\circ (V_s\circ U_s)\circ  U_{s-1}\circ \dots \circ U_{1}=\\
U_s\circ (U_{s-1}\circ V_{s-1})\circ  U_{s-1}\circ \dots \circ U_{1}=U_s\circ U_{s-1}\circ (V_{s-1}\circ  U_{s-1})\circ U_{s-2}\circ \dots \circ U_{1}= \\
= \dots = (U_s\circ U_{s-1}\circ \dots \circ U_{1})\circ B,
\end{gathered}
$$
and 
$$
\begin{gathered}
B\circ  (V_{1}\circ \dots \circ V_{s-1}\circ V_s)=V_1\circ (U_1\circ V_1)\circ V_{2}\circ \dots \circ V_{s-1}\circ  V_{s}=\\
V_1\circ (V_{2}\circ U_{2})\circ  V_{2}\circ \dots \circ V_{s-1}\circ  V_{s}=V_1\circ V_{2}\circ (U_{2}\circ  V_{2})\circ \dots \circ V_{s-1}\circ  V_{s}= \\
= \dots = (V_{1}\circ \dots \circ V_{s-1}\circ V_s)\circ B_s.
\end{gathered}
$$

Further, 
$$
B^{\circ s}=(V_1\circ U_1)\circ (V_1\circ U_1)\circ \dots  \circ (V_1\circ U_1)= V_1\circ B_1^{\circ s-1}\circ U_1 =$$
$$=V_1\circ V_2\circ B_2^{\circ s-2}\circ U_2\circ U_1=\dots 
=(V_1\circ V_2\circ \dots \circ V_{s}) \circ (U_{s}\circ \dots \circ U_2\circ U_1)
$$
and 
$$
B_s^{\circ s}=(U_s\circ V_s)\circ (U_s\circ V_s)\circ \dots  \circ (U_s\circ V_s)= U_s\circ B_{s-1}^{\circ s-1}\circ V_s =$$
$$=U_s\circ U_{s-1}\circ B_{s-2}^{\circ s-2}\circ V_{s-1}\circ V_{s}= \dots
=(U_s\circ U_{s-1}\circ \dots \circ U_{1}) \circ (V_{1}\circ \dots \circ V_{s-1}\circ V_s).\eqno{\Box}
$$

\vskip 0.2cm

The notion of equivalence can be extended to endomorphisms of complex tori. Namely, if $B: R\rightarrow R$ is such an endomorphism, and $B=V\circ U$ is a decomposition of $B$ into a composition of holomorphic maps
$U: R\rightarrow {R}'$ and  $V: {R}'\rightarrow R$ between complex tori,  then
the endomorphism  $U\circ V: {R}'\rightarrow {R}'$ is called an elementary transformation of $B$, and endomorphisms $B: R\rightarrow R$ and $A: T\rightarrow T$
between complex tori are called equivalent if there exists a chain of elementary 
transformations between $B$ and $A$. 
Clearly, an analogue  of Lemma \ref{lem1} holds verbatim for any chain of elementary transformations between endomorphisms of complex tori. Abusing the notation,  
below we will use for equivalent endomorphisms of complex tori  the same symbol  $\sim$ as for equivalent rational functions.

\bt \l{los}
Let  $R$ be a compact Riemann surface  of genus zero or one,   
and $A: \C\P^1\rightarrow \C\P^1$, $B: R\rightarrow R$, and $\pi:R\rightarrow \C\P^1$ holomorphic maps of degree at least two 
such that diagram \eqref{xxuuii}
commutes. Then 
$A$ is a generalized Latt\`es map, unless $R=\C\P^1$ and $B\sim A$.
In more details,
there exist a compact Riemann surface $R_0$ of the same genus as $R$ and 
holomorphic maps    
 $\psi:R\rightarrow R_0,$  $\pi_0:R_0\rightarrow \C\P^1$, and $B_0: R_0\rightarrow R_0$
satisfying the following conditions.

\begin{enumerate}
\item $B_0\sim B$ and $\pi=\pi_0\circ \psi$. 

\item 
The diagram 
\be 
\begin{CD} \l{xxuu}
R @>B>> R \\
@VV  \psi V @VV  \psi  V\\ 
R_0 @> B_0 >> R_0\\
@VV \pi_0 V @VV \pi_0 V\\ 
\C\P^1 @>A >> \ \ \C\P^1\,\  
\end{CD}
\ee
commutes.

\item  The map $\pi_0$ has degree at least two, unless $R=\C\P^1$ and $B\sim A$, and 
the collection 
\be \l{pesa} f=\pi_0, \ \ \ p=B_0, \ \ \ g=A, \ \ \ q=\pi_0\ee
 is a good solution of \eqref{m}.

\item   
The maps $A:\f O_2^{\pi_0}\rightarrow \f O_2^{\pi_0}$ and $B_0:\f O_1^{\pi_0}\rightarrow \f O_1^{\pi_0}$
are minimal holomorphic maps between orbifolds.

\item   
The map $\psi$ is a compositional right factor of $B^{\circ s}$ and a compositional left factor of $B_0^{\circ s}$ 
for some  $s\geq 1.$

\end{enumerate}

\et
\pr 
If the collection \be \l{ur} f=\pi, \ \ \ p=B, \ \ \ g=A, \ \ \ q=\pi\ee is a good  solution of \eqref{m}, we can set  
$$R_0=R, \ \ \ \ \ B_0=B,  \ \ \ \ \ \pi_0=\pi, \ \ \ \ \ \psi=z.$$ Then $A:\f O_2^{\pi_0}\rightarrow \f O_2^{\pi_0}$ and $B_0:\f O_1^{\pi_0}\rightarrow \f O_1^{\pi_0}$
are minimal holomorphic maps 
by Theorem \ref{t1}. The other conditions hold trivially.

Assume now that \eqref{ur} is not a good solution of  \eqref{m}.
Since for solution \eqref{ur}
the third condition of Lemma \ref{good} is always satisfied, this implies that 
$\pi$ and $B$ have a non-trivial common  compositional right factor, that is  there exist 
 a Riemann surface $R^{\prime}$ 
and  holomorphic maps   $$U_1:R\rightarrow R^{\prime},\  \ \ \pi^{\prime}:R^{\prime}\rightarrow \C\P^1, \ \ \ V_1:R^{\prime}\rightarrow R,$$
such that   
\be \l{eli} \pi=\pi^{\prime}\circ U_1 , \ \ \ B=V_1\circ U_1,\ee
and $\deg U_1 \geq 2.$ Furthermore, since $B:R\rightarrow R$ is decomposed as 
$$R\overset{U_1}{\longrightarrow}R^{\prime}\overset{V_1}{\longrightarrow}R,$$
the equality $g(R^{\prime})=g(R)$ holds.

Substituting  \eqref{eli} in the equality $$A\circ \pi=\pi \circ B,$$ we obtain the 
equality 
$$A\circ \pi'= \pi'\circ U_1\circ V_1 $$ and   the commutative diagram 
$$ 
\begin{CD} 
R @>B>> R \\
@VV  U_1 V @VV  U_1  V\\ 
R^{\prime} @> U_1\circ V_1 >> R^{\prime}\\
@VV \pi^{\prime} V @VV \pi^{\prime} V\\ 
\C\P^1 @>A >> \ \ \C\P^1\,.  
\end{CD}
$$
If the solution  $$ f=\pi^{\prime}, \ \ \ p=U_1\circ V_1, \ \ \ g=A, \ \ \ q=\pi^{\prime}$$  of \eqref{m} is still not  good,  
we can perform a similar transformation once again. Since  $\deg U_1\geq 2$ implies that $\deg \pi^{\prime}<\deg \pi,$ it is clear that 
after a finite number of steps we will arrive to diagram \eqref{xxuu}, where $B_0$ is obtained from 
$B$ by a chain of elementary transformations \eqref{seqs} (in the notation of Lemma \ref{lem1}, $B_0=B_s$),
the function $\psi$ has the form
$$\psi=U_{s}\circ \dots \circ U_2\circ U_1,$$ 
and the maps $\pi_0$ and $B_0$ have no non-trivial common  compositional right factor.
 Furthermore, 
$\deg \pi_0=1$ only  if  $R=\C\P^1$ and $B\sim A$.  By Lemma \ref{good}, solution \eqref{pesa} of \eqref{m} is good, and 
applying  Theorem \ref{t1}  
 we obtain that 
$A:\f O_2^{\pi_0}\rightarrow \f O_2^{\pi_0}$ and  $B_0:\f O_1^{\pi_0}\rightarrow \f O_1^{\pi_0}$
are minimal holomorphic maps between orbifolds. Notice that by Proposition \ref{p1} this implies 
that $\chi(\f O_2^{\pi_0})\geq 0$.
Finally, by Lemma  \ref{lem1}, $\psi$ is a compositional factor of $B^{\circ s}$ and a compositional left factor of $B_0^{\circ s}$.
\qed

\vskip 0.2cm

\br\normalfont  
Theorem \ref{los} implies in particular that the problem of describing rational solutions of the functional equation
\be \l{be1} A\circ \pi=\pi\circ B\ee in a sense reduces to the case where  $\chi(\f O_2^{\pi})\geq 0$
 (see \cite{semi} for more details). 
Moreover,  it is shown in the paper \cite{gen}, based on methods of \cite{semi}, that for any good  rational solution of the more general functional equation  
\be \l{be2} A\circ \delta=\pi\circ B,\ee such that $$\deg A\geq 84\, \deg \pi$$ the inequality $\chi(\f O_2^{\pi})\geq 0$ still holds.
The rational functions $\pi$ with $\chi(\f O_2^{\pi})\geq 0$ are characterized by the condition 
 that the genus of the Galois closure of 
$\C(z)/\C(\pi)$ equals zero or one (see \cite{gen}). 
For a detailed description of such functions we refer the reader   
to the paper \cite{gen0}. Notice that  functional equations \eqref{be1} and \eqref{be2} 
naturally arise in arithmetic and dynamics (see e. g. \cite{bilu}, \cite{e}, \cite{ms}, \cite{pj}).
\er

Let us prove now
the chain of implications  $3 \Rightarrow 2 \Rightarrow 1 \Rightarrow 3$ between the conditions of  Theorem  \ref{bas}.
\vskip 0.2cm
\noindent{$3 \Rightarrow 2 $}.
By Proposition \ref{poiu}, for any minimal holomorphic map $A:\, \f O\rightarrow \f O$   between orbifolds 
there exists a holomorphic map $F:\t{\f O}\rightarrow \t{\f O}$ and a homomorphism $\phi:\, \Gamma_{\f O}\rightarrow \Gamma_{\f O}$ such that 
the diagram 
\be \l{diaa}
\begin{CD}
\tt {\f O} @>F>> \tt {\f O}\\
@VV\theta_{\f O}V @VV\theta_{\f O}V\\ 
\f O @>A >> \f O\ 
\end{CD}
\ee
commutes and 
\be \l{hommm} F\circ\sigma=\phi(\sigma)\circ F, \ \ \ \sigma\in \Gamma_{\f O}.\ee
If $\chi(\f O)>0$, then    $\t{\f O}=\C\P^1$ is a compact Riemann surface, so \eqref{xui+} holds for
$$R=\C\P^1, \ \ \ \ B=F, \ \ \ \ \pi=\theta_{\f O}, \ \ \ \ \Gamma=\Gamma_{\f O},$$ and the assumption $\f O\neq\C\P^1$ implies that the  
group $\Gamma$ is non-trivial. Finally, the homomorphism $\phi$ in \eqref{homm2} is an isomorphism by Theorem \ref{las}.

Assume now that $\chi(\f O)=0$ and $\t{\f O}=\C$. Observe first that since  in this case
$A:\, \f O\rightarrow \f O$ is a covering map, 
 the homomorphism $\phi$ in \eqref{homm} is 
a {\it monomorphism}. Indeed,  by Proposition \ref{poiu},  the map $F:\C \rightarrow \C$ is an isomorphism, that is has the form 
\be \l{for} F=az+b, \ \ \ a,b\in \C.\ee   
Thus, $F$ is invertible and hence the 
equality $F\circ \sigma =F$ implies that $\sigma$ is the identity mapping.

Let  now $\Lambda_{\f O}$ be the subgroup of $\Gamma_{\f O}$  generated by translations.
By the classification of groups $\Gamma_{\f O}$ given in the previous section, 
$\theta_{\f O}$ is decomposed as
$$\psi:\C\overset{\psi}{\longrightarrow} \C/\Lambda_{\f O}\cong R\overset{\pi}{\longrightarrow} R/\Gamma\cong \C\P^1,$$ 
where $R=\C/\Lambda_{\f O}$ is a  complex torus and
$\Gamma\cong \Gamma_{\f O}/\Lambda_{\f O}$ is a finite subgroup of $Aut(R).$
Since $\phi$ is a monomorphism, it maps elements of infinite order of $\Gamma_{\f O}$ to elements of infinite order.  
Therefore, $\phi(\Lambda_{\f O})\subset \Lambda_{\f O},$ implying that $F$ descends to a holomorphic map $B: R\rightarrow R$ which makes the 
diagram 
\be 
\begin{CD} \l{gpa3}
\C @>F=ax+b>> \C \\
@VV  \psi V @VV  \psi  V\\ 
R @> B>> R\\
@VV \pi V @VV \pi V\\ 
\C\P^1 @>A >> \ \ \C\P^1
\end{CD}
\ee
commutative. 
 Finally, condition that diagram \eqref{xui+} commutes implies that 
 $B$ commutes with the group $\Gamma$ (see \cite{mil2}, p. 16). Thus,
\eqref{homm2} holds 
for the identical  automorphism $\phi.$ 

\vskip 0.2cm
\noindent{$2 \Rightarrow 1 $}. It is enough to show that if $A,B$ and $\pi$ satisfy the second condition, then $\pi$ is not a compositional right  factor of $B^{\circ s},$ $s\geq 1.$ If $g(R)=1,$ this is obvious, since for any decomposition
$$R\overset{\pi}{\longrightarrow}R^{\prime}\overset{w}{\longrightarrow}R$$
of $B^{\circ s}$, $s\geq 1,$  the genus of $R^{\prime}$ must be equal one. So, assume that 
$R=\C\P^1.$ 

Since  \be \l{ssyl}\pi: \C\P^1\rightarrow \C\P^1/\Gamma\cong \C\P^1\ee
is a Galois covering, 
for any branch point $z_i$, $1\leq i \leq r,$  of $\pi$ there exists a number $d_i$ such that $\pi^{-1}\{z_i\}$ consists of $\vert \Gamma\vert /d_i$ points, and at each of these points the multiplicity of $f$ equals $d_i$.
In other words, 
the orbifold $\f O_1^{\pi}$ is non-ramified.
Since $\C\P^1$ is simply-connected, this implies that $\pi$ is the universal covering of $\f O_2^{\pi}$. Therefore, diagram \eqref{xui+} 
has form \eqref{dia3}, where $\f O=\f O_2^{\pi}$, and Theorem \ref{las} implies that $A:\f O_2^{\pi }\rightarrow\f O_2^{\pi }$ is a minimal
holomorphic map. 
Assume now that 
\be \l{losd} B^{\circ s}=w\circ \pi\ee for some rational function $w$ and $s\geq 1$. Clearly, \eqref{xui+} implies  
\be \l{losd+} A^{\circ s}\circ \pi=\pi \circ B^{\circ s},\ee and substituting \eqref{losd} in \eqref{losd+}, we 
see that 
 \be \l{losd-} A^{\circ s}=\pi\circ w,\ee 
that is $\pi $ is a compositional left factor of $A^{\circ s}$. 
Since $A:\f O_2^{\pi }\rightarrow\f O_2^{\pi }$ is a minimal
holomorphic map, Theorem \ref{serrr} implies that 
$$(A^{\circ s})^*\f O_2^{\pi}=\f O_2^{\pi}.$$ On the other hand, it follows from \eqref{losd-} by Theorem \ref{serrr} that
$$(A^{\circ s})^*\f O_2^{\pi}=(\pi\circ w)^*\f O_2^{\pi}=w^*(\pi^*\f O_2^{\pi})=w^*\f O_1^{\pi}=w^*\C\P^1=\C\P^1.$$
Therefore, $\f O_2^{\pi}=\C\P^1$. However,  for $\Gamma\neq e$ the orbifold $\f O_2^{\pi }$ for quotient map \eqref{ssyl} is ramified.
The contradiction obtained finishes the proof.

\vskip 0.2cm
\noindent{$1 \Rightarrow 3 $}.  
Consider good solution \eqref{pesa} provided by Theorem \ref{los} for the maps $A,B$ and $\pi$ satisfying \eqref{xui}. Observe that $\deg \pi_0\geq 2$ for otherwise the function $\pi$ along with $\psi$ is  a compositional right factor of $B^{\circ s}$ and a compositional left factor of $A^{\circ s}$, 
in contradiction with the assumption.
By Theorem \ref{t1}, $A:\f O_2^{\pi_0}\rightarrow \f O_2^{\pi_0}$ is a minimal holomorphic map, and it follows from $\deg \pi_0\geq 2$ that 
 $\f O_2^{\pi_0}\neq \C\P^1.$  \qed

\br\normalfont The above proof shows that the assumption in the first condition of Theorem \ref{bas}, requiring that $\pi$ is not a compositional right  factor of some iterate of $B$, can be replaced by the assumption  that $\pi$ is not a compositional  left  factor of some iterate of $A$. 

Further, observe that 
for any diagram \eqref{xui+} condition \eqref{homm2}
holds automatically for some {\it homomorphism}  $\phi:\, \Gamma \rightarrow \Gamma $. Moreover, 
if $g(R)=1$, then $\phi$ is an automorphism,
since in this case the commutativity of 
diagram \eqref{xui+} implies that 
 $B$ commutes with  $\Gamma$.
On the other hand, if $g(R)=0$, then, by Theorem \ref{las}, the condition that $\phi$ is an automorphism 
can be replaced by the requirement that $\pi$ and $B$ have no
 common compositional right factor. 

Finally, observe that for surfaces $R$ of genus one  the second condition of Theorem \ref{bas} can be replaced by the condition that there exists a subgroup $\Gamma$ of $ Aut( { \C})$ acting properly discontinuously on $ \C$ whose translation subgroup  is a free group of rank two,
and  a holomorphic map $ F:\, { \C } \rightarrow   { \C}$  
such that diagram \eqref{xui+},
where $\pi:\,  \C\rightarrow  \C/\Gamma$ is the quotient map,  commutes (cf. \cite{mil2}).

\er

\end{section}

\begin{section}{Compositions and decompositions}

For a given orbifold $\f O$  denote by $\f E(\f O)$ the set of rational functions  $A$ such that 
$A:\, \f O\rightarrow \f O$ is a minimal holomorphic map.   In this section we study compositional properties of elements of $\f E(\f O).$

\bt \l{comp} Let $\f O$ be an orbifold and  $U$, $V$  
rational functions of  degree at least two. Assume that $U$ and $V$ are contained in $\f E(\f O)$. Then  the composition $U\circ V$  is also contained in $\f E(\f O).$ In the other direction, if $U\circ V$ is contained in $\f E(\f O)$, 
then  $\nu(U^*\f O)= \nu(\f O)$ and 
 $V:\, \f O\rightarrow U^*\f O$ and $U:\, U^*\f O\rightarrow \f O$ are minimal holomorphic maps.
In particular, whenever
$\nu(\f O)\neq \{2,2,2,2\},$ there exists a M\"obius transformation $\mu$ such that $U\circ \mu$ and $\mu^{-1} \circ V$ are contained in $\f E(\f O).$ 
\et 
\pr 
If $U,V$ are contained in $\f E(\f O)$, then Corollary \ref{serka0} obviously  implies that the composition $U\circ V$  is also contained in $\f E(\f O)$.

In the other direction, assume that $U\circ V\in \f E(\f O)$, and 
set $\f O^{\prime}=U^*\f O.$ 
Since by Corollary \ref{indu2}  \be \l{gopp} U: \f O^{\prime}\rightarrow \f O, \ \ \ \ \ V: \f O\rightarrow \f O^{\prime}\ee
are minimal holomorphic maps between orbifolds, 
 we have: 
\be \l{xer} \nu(\f O)\leq \nu (\f O^{\prime} )\leq \nu(\f O).\ee 
Furthermore, by Proposition \ref{p1}, the inequalities 
$$\chi(\f O)\leq \chi(\f O^{\prime})\deg V, \ \ \ \ \ \chi(\f O^{\prime})\leq \chi(\f O)\deg U$$ hold.
Therefore, 
$$\chi(\f O)\leq \chi(\f O^{\prime})\deg V\leq \chi(\f O)\deg U\deg V,$$ implying that  
$\chi(\f O^{\prime})=0$ whenever $\chi(\f O)=0$, and $\chi(\f O^{\prime})>0$ whenever $\chi(\f O)>0.$

Assume first that  $\chi(\f O)=0$.
Then a direct analysis of the table
\begin{table}[H]
\caption*{Table 1}
\centering
\begin{tabular}{ | l | c | c | c | r | } 
\hline	\ &  \{2,2,2,2\}& \{3,3,3\} & \{2,4,4\} & \{2,3,6\}\\ 
\hline \{2,2,2,2\} & $\leq$ & \ & $\leq$ & $\leq$ \\ 
\hline \{3,3,3\} &  &$\leq$  &  & $\leq$ \\ 
\hline \{2,4,4\} &  &  &  $\leq$&  \\ 
\hline \{2,3,6\} &  &  &  &  $\leq$\\ 
\hline \end{tabular}
\end{table}

\noindent listing all $\nu(\f O_1)$ and $\nu(\f O_2)$ such that  
 $$\chi(\f O_1)=\chi(\f O_2)=0$$ and $\nu(\f O_1)\leq \nu(\f O_2)$,  
shows that \eqref{xer} is possible only if $\nu (\f O^{\prime} )=\nu(\f O)$. 

If $\chi(\f O)>0$ the proof can be done as follows (cf. \cite{semi}, Corollary 5.1). Since maps \eqref{gopp} 
are minimal holomorphic maps, it follows from  Proposition \ref{poiu} that there exist
rational functions $F_U$ and $F_V$ which make the diagram    
\begin{equation*} 
\begin{CD} 
\C\P^1 @>F_V>> \C\P^1 @>F_U>>\C\P^1 \\
@VV\theta_{\f O}  V @VV \theta_{\f O^{\prime}} V @VV \theta_{\f O} V\\ 
\C\P^1 @> V >> \C\P^1 @>U>>\C\P^1 \\
\end{CD}
\end{equation*}
commutative and satisfy 
$$F_V\circ\sigma=\phi_V(\sigma)\circ F_V, \ \ \sigma\in \Gamma_{\f O}, \ \ \ \ 
F_U\circ\sigma=\phi_U(\sigma)\circ F_U, \ \ \ \sigma\in \Gamma_{\f O^{\prime}},$$
for some homomorphisms 
$$\phi_V:   \Gamma_{\f O}\rightarrow \Gamma_{\f O^{\prime}}, \ \ \  \phi_U:   \Gamma_{\f O^{\prime}}\rightarrow \Gamma_{\f O}.$$
Since the function $F_U\circ F_V$ makes the diagram 
\begin{equation*} 
\begin{CD} 
\C\P^1 @>F_U\circ F_V>>\C\P^1 \\
@VV\theta_{\f O}  V @VV  \theta_{\f O} V\\ 
\C\P^1 @> U \circ V>>\C\P^1 \\
\end{CD}
\end{equation*}
 commutative, 
 Theorem \ref{las} implies that the composition of homomorphisms $$\phi_{U}\circ \phi_{V}: \Gamma_{\f O}\rightarrow \Gamma_{\f O}$$ is an 
automorphism. Therefore, $\Gamma_{\f O^{\prime}}\cong\Gamma_{\f O}$, implying that  $\nu(\f O^{\prime})=\nu(\f O)$.

Finally, if  $\nu(\f O)\neq \{2,2,2,2\}$, the orbifolds $\f O$ and $\f O^{\prime}$ have at most three singular points, implying that we can find 
$\mu$ as required.  
\qed
\vskip 0.2cm
In a sense, Theorem \ref{comp} 
reduces the study of generalized Latt\`es maps  to the study of indecomposable maps. Recall that a rational function $A$  is called {\it indecomposable} if the equality $A=U\circ V$, where $U$ and $V$ are rational functions,  implies that at least one of the functions $U$ and $V$ has degree one.
Clearly, any rational function   $A$ of degree at least two can be decomposed into a composition 
\be \l{max} A=A_1\circ A_2\circ \dots \circ A_l\ee
of indecomposable rational functions of degree at least two. Such decompositions are called {\it maximal}.

\bc Let $\f O$ be an orbifold  whose signature is distinct from $\{2,2,2,2\}$. Then 
any rational function $A$ of degree at least two contained in $\f E(\f O)$ has a maximal decomposition 
whose elements are contained in 
$\f E(\f O)$.\qed
\ec
\pr Indeed, if $A$ is indecomposable we have nothing to prove. Otherwise, $A=U\circ V$ for some rational functions $U$ and $V$, and 
changing  $U$ to $U\circ \mu$ and $V$ to $\mu^{-1} \circ V$,
where $\mu$ is a M\"obius transformation provided by Theorem \ref{comp}, without loss of generality we may assume that $U,V\in  \f E(\f O)$.
Continuing in this way we will obtain the required maximal decomposition. \qed

\bc \l{sosis} Let $\f O$ be an orbifold  whose signature is distinct from $\{2,2,2,2\}$. Assume that $A\in \f E(\f O)$ and $B\sim A$. Then $B$ is conjugate to some $B'\in \f E(\f O)$.
\ec
\pr By Theorem \ref{comp},  the statement is true for any elementary transformation of $A$.
 It follows now from the definition of the equivalence $\sim$ that  it is true for any $B\sim A$. \qed

\vskip0.2cm


\bc \l{korova2} Let $A$ be a Latt\`es map and $B\sim A$. Then $B$ is a Latt\`es map. \ec
\pr It follows from Theorem \ref{comp} and Corollary \ref{serka0} that if $A=U\circ V$  is contained in $\f E(\f O)$, then the elementary transformation
 $V\circ U$  is contained in $\f E(U^*\f O)$. Moreover, since $\nu(U^*\f O)=\nu(\f O),$ if $\chi(\f O)=0,$
then $\chi(U^*\f O)=0.$ Therefore, 
if $A=U\circ V$ is a Latt\`es map, then  $V\circ U$ is also a Latt\`es map.
 \qed

\vskip0.2cm

For orbifolds $\f O_1, \f O_2, \dots ,\f O_s$ define the orbifold 
$\f {O}=\LCM(\f O_1,\f O_2, \dots ,\f O_s)$ by the condition 
$$ \nu(z)=\LCM\Big(\nu_1(z),\nu_2(z), \dots, \nu_s(z)\Big), \ \ \ z\in \C\P^1.$$

\bt \l{krek}  Let $\f O_1,$ $\f O_2,\dots ,\f O_s$ and  $\f O_1^{\prime},$ $\f O_2^{\prime},\dots ,\f O_s^{\prime}$ 
be orbifolds, and $A$ a rational function such that the maps $A:\, \f O_i\rightarrow \f O_i^{\prime}$, $1\leq i \leq s,$ 
are holomorphic maps (resp. minimal holomorphic maps, covering maps) between orbifolds. Then $$A:\,  \LCM(\f O_1,\f O_2, \dots ,\f O_s)\rightarrow \LCM(\f O_1^{\prime},\f O_2^{\prime}, \dots ,\f O_s^{\prime})$$  is also a holomorphic map (resp. a minimal holomorphic map, a covering map)  between orbifolds.
\et
\pr 
In order to prove the first part of the proposition, it is enough to observe that the conditions 
$$\nu_{i}^{\prime}(A(z)) \mid \nu_{i}(z)\deg_zA, \ \ \ 1\leq i \leq s,$$ 
imply the condition 
$$\LCM\Big(\nu_{1}^{\prime}(A(z)),\nu_{2}^{\prime}(A(z)), \dots , \nu_{s}^{\prime}(A(z))\Big)\mid $$ 
$$ \LCM\Big(\nu_{1}(z)\deg_zA,\, \nu_{2}(z)\deg_zA,\, \dots , \nu_{s}(z)\deg_zA\Big)=$$ 
$$\LCM\Big(\nu_{1}(z),\nu_{2}(z), \dots , \nu_{s}(z)\Big)\deg_zA.  $$

In order to prove the second part,  we must show that if 
$$\nu_{i}^{\prime}(A(z)) = \nu_{i}(z)\GCD\Big(\nu_{i}^{\prime}(A(z)),\deg_zA\Big), \ \ \ 1\leq i \leq s,$$
then  
 $$\LCM\Big(\nu_{1}^{\prime}(A(z)),\nu_{2}^{\prime}(A(z)), \dots , \nu_{s}^{\prime}(A(z))\Big)=\LCM\Big(\nu_{1}(z),\nu_{2}(z), \dots , \nu_{s}(z)\Big)\times $$
\be \l{a}\times\GCD\bigg(\LCM\Big(\nu_{1}^{\prime}(A(z)),\nu_{2}^{\prime}(A(z)), \dots , \nu_{s}^{\prime}(A(z))\Big),\deg_zA\bigg).\ee
Let $p$ be an arbitrary prime number and $z\in \C\P^1$. Set 
$$b_i=\ord_p\nu_{i}^{\prime}(A(z)), \ \ \ a_i=\ord_p\nu_{i}(z), \ \ \  c=\ord_p\deg_zA, \ \ \ 1\leq i \leq s.$$
Considering the orders at $p$ of the numbers in the left and the right sides of equality \eqref{a}, we see that we must prove the following statement: if $a_i,$ $b_i$, $1\leq i \leq s,$ and $c$ are integer non-negative numbers  such that 
\be \l{a1} 
b_i=a_i+\min\{c,b_i\}, \ \ \  1\leq i \leq s,
\ee
then 
\be \l{a2} 
\max\limits_{i}\{b_i\}=\max\limits_{i}\{a_i\}+\min\{c,\max\limits_{i}\{b_i\}\}.
\ee

Let $I_1$ (resp. $I_2$) be the subset of $\{1,2,\dots s\}$ consisting of indices $i$ such that  $c\leq b_i$ (resp.  $c>b_i$). 
Clearly, we have:
$$ \max\limits_{i}\{b_i\}=\max\Big\{\max\limits_{i\in I_1}\{b_i\}, \max\limits_{i\in I_2}\{b_i\}\Big\}.$$
For each $i,$ $1\leq i \leq s,$ equality \eqref{a1} implies that $b_i=a_i+c,$ if $i\in I_1$, and $a_i=0$, if $i\in I_2.$
If $c> \max\limits_{i}\{b_i\}$, that is the set $I_1$ is empty, then $\max\limits_{i}\{a_i\}=0,$ and hence \eqref{a2} holds. 
On the other hand,
if  $c\leq \max\limits_{i}\{b_i\}$, then $I_1$ is non-empty and for an arbitrary $i_0\in I_1$  we have $b_{i_0}=a_{i_0}+c$, implying that for any $i\in I_2$ the inequality
$$b_i<c\leq c+a_{i_0}=b_{i_0}\leq \max\limits_{i\in I_1}\{b_i\}$$ holds.
Thus, 
$$\max\limits_{i\in I_2}\{b_i\}< 
\max\limits_{i\in I_1}\{b_i\}$$
and hence 
$$ \max\limits_{i}\{b_i\}=\max\limits_{i\in I_1}\{b_i\}=\max\limits_{i\in I_1}\{a_i+c\}=\max\limits_{i\in I_1}\{a_i\}+c.$$
Furthermore, since $a_i=0$ whenever $i\in I_2,$ we have:
$$\max\limits_{i\in I_1}\{a_i\}=\max\limits_{i}\{a_i\}.$$
Therefore, if $c\leq \max\limits_{i}\{b_i\}$, then 
$$ \max\limits_{i}\{b_i\}=\max\limits_{i}\{a_i\}+c,$$
as required.

Finally, since a   minimal holomorphic map $f:\, \f O\rightarrow \f O^{\prime}$ is a covering map if and only if $\deg_z A\vert \nu^{\prime}(A(z))$ for any $z\in \C\P^1$,  
in order to prove the last part of the theorem it is enough to observe that the conditions $$\deg_z A\mid \nu^{\prime}_i(A(z)),   \ \ \ \ \ 1\leq i \leq s,\ \ \ \ \ z\in \C\P^1,$$ imply the condition  
$$\deg_z A\mid \LCM\Big(\nu_{1}^{\prime}(A(z)),\nu_{2}^{\prime}(A(z)), \dots , \nu_{s}^{\prime}(A(z))\Big), \ \ \ \ z\in \C\P^1. \eqno{\Box}$$

\bc \l{lat} Let $A$ be a rational function of degree at least two, and $\f O$ an orbifold such that the function $A^{\circ l}$ is contained in $\f E(\f O)$ for some $l\geq 2.$ Then, 
unless the signature of $\f O$ is $\{2,2\},$ $\{3,3\}$, $\{2,2, 2\},$ or $\{2,2, 4\},$ the function $A$ is also contained in $\f E(\f O)$. 
\ec
\pr Set $\f O'=A^*\f O.$  Applying 
Theorem \ref{comp} to the decomposition $$A^{\circ l}= A\circ A^{\circ (l-1)}$$ we see that $\nu(\f O')=\nu(\f O)$ and 
the maps \be\l{goper}  A:\f O'\rightarrow \f O, \ \ \ \ \ A^{\circ (l-1)}:\f O\rightarrow \f O'\ee  are minimal holomorphic maps. 
In particular, in order to show that $A\in \f E(\f O)$ it is enough to prove that $\f O'=\f O.$
Since \eqref{goper} are minimal holomorphic maps, applying Corollary \ref{serka0} to the decomposition 
$$A^{\circ l}= A^{\circ (l-1)}\circ A,$$ we see that 
$A^{\circ l}\in \f E(\f O').$
It follows now from  Theorem \ref{krek} that $A^{\circ l}\in \f E(\t{\f O}),$
where $\t{\f O}=\LCM(\f O,\f O').$ 
 However, this implies that $\chi(\t{\f O})\geq 0$, and it is easy to see that if $\f O$ and $\f O'$ are two orbifolds of non-negative Euler characteristic such that $\nu(\f O')=\nu(\f O)$ and  $\chi(\t{\f O})\geq 0$, then $\f O'=\f O$, unless the signature of 
  $\f O$ is  $\{2,2\},$ $\{3,3\}$,  $\{2,2, 2\},$ or $\{2,2, 4\}.$ 
  Indeed, assume that, say, $\nu(\f O)=\{2,2,n\},$ $n\geq 2.$ 
Since $\chi(\t{\f O})\geq 0$, if $c(\f O')\neq 	c(\f O)$, then $c(\t{\f O})$ contains four points  and $\nu(\t{\f O})=\{2,2,2, 2\},$ so that $n=2.$ On the other hand, 
if $c(\f O')= 	c(\f O)$ but $\f O'\neq \f O$,  then $\nu(\t{\f O})=\{2,d,d\},$ where $d=\LCM(2,n),$ implying that $n=4.$ Other signatures can be considered similarly.
	\qed

\vskip 0.2cm

Notice that Corollary \ref{lat} implies in particular the following statement.

\bc \l{korova1} Let $A$ be a rational function of degree at least two such that some iterate $A^{\circ l}$ is a Latt\`es map. Then $A$ is a Latt\`es map. \qed
\ec

\end{section}

\begin{section}{Generalized Latt\`es maps for the signatures $\{n,n\}$ and $\{2,2,n\}$}
In this section we describe  minimal holomorphic maps   $A:\, \f O\rightarrow \f O$ for orbifolds $\f O$ 
with signatures $\{n,n\}$ and $\{2,2,n\}$. To be definite, we normalize considered orbifolds   by the conditions \be \l{norm} \nu(0)=n, \ \ \ \nu(\infty)=n, \ \ \ n\geq 2,\ee 
and 
 \be \l{norm2}\nu(-1)=2, \ \ \ \nu(1)=2, \ \ \ \nu(\infty)=n, \ \ \  n>2.\ee
For the orbifold $\f O$ defined by \eqref{norm} 
 the corresponding group 
$\Gamma_{\f O}$ is a cyclic group $C_n$  generated by 
 \be \l{zach} \alpha: z\rightarrow e^{2\pi i/n} z,\ee 
 and \be \l{xria} \theta_{\f O}=z^n.\ee 
For  $\f O$ defined by \eqref{norm2} 
the group 
$\Gamma_{\f O}$ is a dihedral group $D_n$ generated by 
\be \l{59} \alpha: z\rightarrow e^{2\pi i/n} z, \ \ \  \ \ \ \beta: z\rightarrow \frac{1}{z},\ee and 
\be \l{xriakk}  \theta_{\f O}=\frac{1}{2}\left(z^n+\frac{1}{z^n}\right).\ee
Notice that the assumption $n>2$ in \eqref{norm2} is due to the fact that the description of the group $Aut(D_{2n})$ in the case $n=2$ is different from the general case.  The case $n=2$ can be analyzed by the method of the seventh section.

By Theorem \ref{las}, $A:\, \f O\rightarrow \f O$ is a 
minimal holomorphic map   for an orbifold $\f O$ with $\chi(\f O)>0$  if and only if the solution of \eqref{m} provided by 
the  commutative diagram 
\be \l{di}
\begin{CD}
\C\P^1 @>F>> \C\P^1\\
@VV\theta_{\f O}V @VV\theta_{\f O}V\\ 
\C\P^1 @>A >> \C\P^1
\end{CD}
\ee
is good, or, equivalently,  the homomorphism   $\phi:\, \Gamma_{\f O}\rightarrow \Gamma_{\f O}$ 
defined by  
\be \l{di+} F\circ\sigma=\phi(\sigma)\circ F, \ \ \ \sigma\in \Gamma_{\f O},\ee 
is an automorphism.
Thus, the problem of describing minimal holomorphic map   $A:\, \f O\rightarrow \f O$ for orbifolds $\f O$ 
defined by \eqref{norm} and \eqref{norm2}  essentially is equivalent to the problem of describing 
good solutions 
of the functional equations 
\be \l{xori1} 
A\circ z^n=z^n\circ F, \ \ \ n\geq 2,
\ee
and 
\be\l{xori2} 
A\circ \frac{1}{2}\left(z^n+\frac{1}{z^n}\right)=\frac{1}{2}\left(z^n+\frac{1}{z^n}\right)\circ F,\ \ \ n>2,
\ee
or, equivalently, to the problems of describing
$F$ satisfying \eqref{di+} for automorphisms $\phi$ of $\Gamma_{\f O}=C_n$ and  $\Gamma_{\f O}=D_{2n}.$ 

Abusing the notation, we will say that a couple of rational functions $A$, $F$ 
is  a good solution of \eqref{xori1} if the functions $A$, $z^n,$ $z^n,$ $F$ form a good solution of \eqref{m}.
A good solution of \eqref{xori2} is defined similarly.

\bt \l{ori1} A couple of rational functions $A$, $F$ 
is  a good solution of \eqref{xori1} if and only if 
$A=z^rR^n(z)$ and $F=z^rR(z^n)$, where $R\in \C(z)$ and 
$\GCD(r,n)=1$.  In particular, any minimal holomorphic map   $A:\, \f O\rightarrow \f O$ for $\f O$ defined by \eqref{norm} has the above form.
\et
\pr
Since for $\Gamma_{\f O}$ generated by \eqref{zach} any
automorphism $\phi: \Gamma_{\f O}\rightarrow  \Gamma_{\f O}$ has the form \be \l{trann} \phi(\alpha)=\alpha^{\circ r},\ \ \ 1\leq r \leq n-1, \ \ \ \GCD(n,r)=1,\ee
a rational function $F$ satisfies \eqref{di+} if and only if for some $r$ coprime with $n$ the function $F/z^r$ is $\Gamma_{\f O}$-invariant, that is $F/z^r$ 
is a rational function in $z^n$. Thus, $F$ satisfies \eqref{di+} if and only if 
$F=z^rR(z^n)$, where $R\in \C(z)$ and $\GCD(r,n)=1.$ 
Finally, it follows from
$$A\circ z^n=z^n\circ z^rR(z^n)=z^rR^n(z)\circ z^n$$
that $A$ makes diagram \eqref{di} commutative if and only  if $A=z^rR^n(z)$. \qed


Notice that $A=z^rR^n(z)$ and $F=z^rR(z^n)$ make diagram \eqref{di} commutative for any $r\geq 0$,
not necessarily coprime with $n$. However, if $\GCD(r,n)>1,$ the homomorphism $\phi$ has a non-trivial kernel, and  $A:\, \f O\rightarrow \f O$ is a holomorphic map but not a minimal holomorphic map.

\bc \l{esa} Let $A$, $F$ be a good solution of \eqref{xori1} and $m=\deg A=\deg F.$ Then $m\geq n$, unless 
$F=cz^{\pm m}$  and $A=c^nz^{\pm m}$, where $c\in \C.$
\ec
\pr Indeed, if a rational function $R$ has a zero or a pole  distinct from $ 0$ and $\infty$, then the degree of the function $F=z^rR(z^n)$ is at least $n.$ Otherwise, $F=cz^{\pm m}$ implying that $A=c^nz^{\pm m}$.
\qed

\bc \l{l2} Let $A$ be a rational function of degree $m\geq 2$  such that \linebreak
$A:\f O\rightarrow \f O$ is a minimal holomorphic map between orbifolds with
$\nu(\f O)=\{n,n\}$, $n\geq 2.$ 
Then $m\geq n$,  unless $A$ is conjugate to $z^{\pm m}$. \qed
\ec 


Denote by $\mathfrak T$ the set of rational functions commuting with the involution
$$\beta: z\rightarrow \frac{1}{z}.$$ 
Since the equality $G(z)G(1/z)=1,$ where $G$ is a rational function, 
implies that   $a\in \C\P^1$ is a zero of $G$ of order $k$ if and only if $1/a$ is a pole of $G$ of order $k$, it is easy to see that elements of $\mathfrak T$
have the form
$$G=\pm z^{\pm l_0}\frac{(z-a_1)^{l_1}(z-a_2)^{l_2}\dots (z-a_s)^{l_s}}{(a_1z-1)^{l_1}(a_2z-1)^{l_2}\dots (a_sz-1)^{l_s}},$$ where $a_1,a_2,\dots a_s\in \C\setminus\{0\}$ and $l_0, l_1,l_2,\dots l_s\in \N.$ 

\bt \l{zxcv2}  A couple of rational functions $A$, $F$ 
is  a good solution of \eqref{xori2} if and only if 
$ F=\v z^rR(z^n)$ and 
\be \l{uri} A=\frac{\v^n}{2}\left(z^rR^n(z)\circ (z+\sqrt{z^2-1})+z^rR^n(z)\circ (z-\sqrt{z^2-1})\right),
\ee 
where $R\in \mathfrak T$,   
$\GCD(r,n)=1$, and $\v^{2n}=1.$   
 In particular, any minimal holomorphic map   $A:\, \f O\rightarrow \f O$ for $\f O$ defined by \eqref{norm2} has the above form.
\et

\pr Since an automorphism $\phi$ of the group $\Gamma_{\f O}$ generated by \eqref{59} 
maps any element of order $n$ of the group $\Gamma_{\f O}=D_{2n}$ to an element of order $n$, and $n> 2$, 
equality \eqref{trann} still holds.
On the other hand, since $\phi$ maps $\beta$  
to an element of order two not belonging to the subgroup generated by $\alpha$, we have: \be \l{trann2} \phi(\beta)=\alpha^{\circ k}\circ\beta=e^{2\pi ik/n}z \circ  \frac{1}{z}, \ \ \ 0 \leq k \leq n-1.\ee
It was shown above that condition  \eqref{trann} holds if and only if $F= z^rR(z^n)$, where   $R\in \C(z)$ and $\GCD(r,n)=1$.
On the other hand,  condition \eqref{trann2} holds if and only if
\be \l{ebus} F(1/z)=e^{\frac{2\pi i}{n}k}\frac{1}{F(z)},\ee 
or equivalently if and only if $e^{-\frac{\pi i}{n}k}F\in \mathfrak T,$ $0 \leq k \leq n-1.$
This implies that $F$ satisfies \eqref{di+} for some automorphism $\phi$ of $\Gamma_{\f O}$ if and only if 
\be \l{larek} F=\v z^rR(z^n),\ee where $R\in \mathfrak T,$ $\v^{2n}=1,$   and $\GCD(r,n)=1$.

Finally, if 
\be \l{xom1} A\circ \frac{1}{2}\left(z^n+\frac{1}{z^n}\right)
=\frac{1}{2}\left(z^n+\frac{1}{z^n}\right)\circ \v z^rR(z^n),\ee then 
it follows from 
$$A\circ \frac{1}{2}\left(z^n+\frac{1}{z^n}\right)=A\circ \frac{1}{2}\left(z+\frac{1}{z}\right)\circ z^n$$ and 
$$\frac{1}{2}\left(z^n+\frac{1}{z^n}\right)\circ \v z^rR(z^n)=\frac{\v^n}{2}\left(z+\frac{1}{z}\right)\circ  z^rR^n(z)\circ z^n,$$
that
\begin{multline} \l{xom2} A\circ \frac{1}{2}\left(z+\frac{1}{z}\right)=\frac{\v^n}{2}\left(z+\frac{1}{z}\right)\circ z^rR^n(z)=\\ \frac{\v^n}{2}\left(z^rR^n(z)+z^rR^n(z)\circ \frac{1}{z}\right).\end{multline} 
Substituting now $z$ by $ z+\sqrt{z^2-1}$ in the left and the right sides of the last equality 
 we obtain \eqref{uri}. On the other hand, if \eqref{uri} holds, then 
substituting $z$ by $$\frac{1}{2}\left(z+\frac{1}{z}\right)$$ we obtain \eqref{xom2} and \eqref{xom1}. \qed

\bc \l{esa2} Let $A$, $F$ be a good solution of \eqref{xori2} and $m=\deg A=\deg F.$ Then $m\geq n+1$, unless 
$F=\v z^{\pm m},$ where  $\v^{2n}=1,$  and $A= \v^nT_{m}$.
\ec
\pr Indeed, if a rational function $R\in \mathfrak T$  has say a zero $a$ distinct from $ 0$ and $\infty$, then it has a pole $1/a$ also distinct from  $ 0$ and $\infty$. Therefore, the function $F=\v z^rR(z^n)$ has the degree at least $n+r\geq n+1.$

On the other hand, if $R\in \mathfrak T$ has  no zeroes or poles distinct from  
$0$ and $\infty$, then $R=\pm z^{\pm l}$, $l\geq 1.$ Therefore,
$F=\v z^{\pm m}$,  where  $\v^{2n}=1,$ and 
the well known identity 
$$T_m(z)=\frac{1}{2}\left((z+\sqrt{z^2-1})^m+(z-\sqrt{z^2-1})^m\right)$$ implies that $A=\v^nT_{m}$. \qed

\bc \l{l2+} Let $A$ be a rational function of degree $m\geq 2$  such that \linebreak
$A:\f O\rightarrow \f O$ is a minimal holomorphic map between orbifolds with $\nu(\f O)=\{2,2,n\}$, $n> 2.$
Then $m\geq n+1$,  unless $A$ is conjugate to $\pm T_m$. \qed
\ec

In conclusion of this section, we provide 
 a description of good solutions of the equation 
\be\l{xori3} 
A\circ T_n=T_n\circ B,\ \ \ n> 2,
\ee
based on
Theorem \ref{zxcv2}.

\bt  A couple of rational functions $A$, $B$ 
is  a good solution of \eqref{xori3} if and only if 
\be \l{bur} B=\frac{1}{2}\left(z^rR(z^n)\circ (z+\sqrt{z^2-1})+z^rR(z^n)\circ (z-\sqrt{z^2-1})\right),
\ee and 
\be \l{ira} A=\frac{1}{2}\left(z^rR^n(z)\circ (z+\sqrt{z^2-1})+z^rR^n(z)\circ (z-\sqrt{z^2-1})\right),
\ee 
where $R\in \mathfrak T$ and  
$\GCD(r,n)=1$.   
\et
\pr Assume that $A$, $B$ 
is  a good solution of \eqref{xori3}.
Observe that for $n>2$ the orbifold $\t{\f O}=\f O_1^{T_n}$ is defined  
by the equalities 
$$\t\nu(-1)=2,\ \ \ \ \ \t\nu(1)=2.$$ 
Since   $B:\t{\f O} \rightarrow \t{\f O} $ is a minimal holomorphic  map
by Theorem \ref{t1}, this implies by Proposition \ref{poiu} that we can complete 
\eqref{xori3} 
to the diagram
\be \l{olm} 
\begin{CD} 
\C\P^1  @>F>> \C\P^1  \\
@VV  \frac{1}{2}\left(z+\frac{1}{z}\right) V @VV  \frac{1}{2}\left(z+\frac{1}{z}\right) V\\ 
\C\P^1 @> B>> \C\P^1\\
@VV T_n V @VV T_n V\\ 
\C\P^1 @>A >> \ \ \C\P^1\,.
\end{CD}
\ee
Furthermore, since $A:\f O_2^{T_n} \rightarrow \f O_2^{T_n}$ is also a minimal holomorphic  map, and $\f O_2^{T_n}$ coincides with $\f O$ defined by \eqref{norm2}, the solution $A,$ $F$ of \eqref{xori2} 
induced by \eqref{olm} is good by Theorem \ref{las}, so that   
equalities \eqref{uri} and \eqref{larek} hold. 

Applying Proposition \ref{poiu} to the upper square of diagram \eqref{olm},
we see that $F$ maps the subgroup generated by $\beta$ to itself. Thus, $k=0$ in \eqref{ebus}, and hence  
$F=z^rR(z^n)$, implying that \eqref{uri}  reduces to \eqref{ira}.
Moreover, 
substituting  $z$ by $ z+\sqrt{z^2-1}$ in the left and the right sides of the equality 
$$ B\circ \frac{1}{2}\left(z+\frac{1}{z}\right)=\frac{1}{2}\left(z+\frac{1}{z}\right)\circ  z^rR(z^n)=\frac{1}{2}\left( z^rR(z^n)+z^rR(z^n)\circ \frac{1}{z}\right),$$ 
we obtain \eqref{bur}.

In the other direction, assume that $A$ and $B$ are given by \eqref{bur} and \eqref{ira}. Since 
in this case the function $F=z^rR(z^n)$ satisfies the equalities \eqref{xori2} and
$$B\circ \frac{1}{2}\left(z+\frac{1}{z}\right)=\frac{1}{2}\left(z+\frac{1}{z}\right)\circ  F,$$
the functions $A$ and $B$ satisfy  \eqref{xori3}. Finally, Lemma \ref{good} implies that $A$, $B$ is a good solution 
 of \eqref{xori3}. Indeed, since $A$, $F$ is a good solution of \eqref{xori2}, the curve 
$$A(x)- \left(y^n+\frac{1}{y^n}\right)=0$$ is irreducible. Therefore, since 
$T_n(z)$ is a compositional left factor of $\left(z^n+\frac{1}{z^n}\right),$ the curve 
$$A(x)- T_n(y)=0$$ is also irreducible, implying that  $A$, $B$ is a good solution 
 of \eqref{xori3}. 
 \qed

\bc \l{esa3} Let $A$, $B$ be a good solution of \eqref{xori3} and $m=\deg A=\deg B.$ Then $m\geq n+1$, unless 
$B=\pm T_{m}$  and $A=(\pm 1)^nT_{m}$.
\ec
\pr Since a good solution $A$, $B$  of \eqref{xori3} induces a good  solution $A,$ $F$ of \eqref{xori2}, 
the corollary is obtained by  a modification of the proof of Corollary \ref{esa2}, taking into account that 
$k=0$ in \eqref{ebus}. \qed

\end{section}

\begin{section}{Orbifold $\f O_0^A$}
Let $A$ be a rational function of degree at least two. In this section we study the totality of orbifolds
$\f O$ such that   $A:\,  \f O\rightarrow \f O$ is a  minimal holomorphic  map, and prove Theorem \ref{uni}.

If $A$ is an ordinary Latt\`es maps, then an orbifold $\f O$ such that $A:\,  \f O\rightarrow \f O$ is a  covering map,  is 
defined in a unique way by dynamical properties of $A$ (see \cite{mil2}). We start by reproving the uniqueness of $\f O$   
using Theorem \ref{krek}. 

\bt Let be a rational function $A$ of degree at least two. Then there exists at most one orbifold $\f O$ of zero Euler characteristic such that $A:\f O\rightarrow \f O$ is a  minimal holomorphic  map between orbifolds. 
\et
\pr 
Assume that $\f O_1,$ $\f O_2$ are two such orbifolds, and set $\f O=\LCM(\f O_1,\f O_2).$  By Proposi\-tion \ref{p1},  $A:\f O_1\rightarrow \f O_1$ 
and $A:\f O_2\rightarrow \f O_2$ are covering maps between orbifolds.  
Therefore, $A:\f O\rightarrow \f O$ is also a covering map, by Theorem \ref{krek}. Thus, $\chi(\f O)=0$. However,  it is easy to see that whenever  $\nu(\f O_1)$ and $\nu(\f O_2)$ belong to  list \eqref{list} 
 the equality $\chi(\f O)= 0$ implies the equality $\f O_1=\f O_2.$ \qed 

\vskip 0.2cm

In general, there might be more than one orbifold $\f O$ such that $A:\, \f O \rightarrow \f O$ is a minimal holomorphic map  between orbifolds, and even infinitely many such orbifolds.  The last phenomenon occurs for the functions $z^{\pm d}$ and $\pm T_{d}$, which play a  special role in the theory.
Namely,  $z^{\pm d}:\f O\rightarrow \f O$ is a minimal holomorphic map for any $\f O$ defined by the conditions \be \l{raz} 
\nu(0)=\nu(\infty)=n, \ \ \ \ n\geq 2, \ \ \ \  \GCD(d,n)=1,\ee and  $\pm T_{d}:\f O\rightarrow \f O$ is a minimal holomorphic map  for any $\f O$ defined by  the conditions \be \l{dva} 
\nu(-1)=\nu(1)=2, \ \ \ \nu(\infty)=n,\ \ \ \ n\geq 1, \ \ \ \  \GCD(d,n)=1.\ee 
Indeed, it is enough to check condition \eqref{uu}  only at points of the finite set 
\be \l{krysa} c(\f O)\cup A^{-1}(c(\f O)),\ee since at other points it  holds trivially, 
and at  points of \eqref{krysa} this condition  holds by the well-known ramification properties of $z^{\pm d}$ and $\pm T_{d}$.

Notice that for odd $d$, additionally, $\pm T_{d}:\f O\rightarrow \f O$ is a minimal holomorphic map 
for $\f O$ defined 
by 
\be \l{tri} 
\nu(1)=2, \ \ \ \nu(\infty)=2,\ee 
or \be \l{chet} 
\nu(-1)=2, \ \ \ \nu(\infty)=2.\ee   

\bt \l{214} Let $\f O$ be an orbifold distinct from the non-ramified sphere. 
\begin{enumerate}

\item The map  $z^{\pm d}:\f O\rightarrow \f O$, $d\geq 2,$ is a minimal holomorphic map between orbifolds if and only if $\f O$ is defined by conditions \eqref{raz}. 

\item The map $\pm T_{d}:\f O\rightarrow \f O$, $d\geq 2,$ is a minimal holomorphic map between orbifolds 
if and only if either $\f O$ is defined 
by conditions \eqref{dva},  or $d$ is odd and 
 $\f O$ is defined 
by conditions
\eqref{tri} or  
\eqref{chet}.

\end{enumerate}

\et 
\pr  We prove the theorem for $\pm T_d$. For $z^{\pm d}$ the proof is similar. Assume that  $\pm T_{d}:\f O\rightarrow \f O$ is a minimal holomorphic map between orbifolds, and set $\f O_n$ equal LCM of the orbifolds  $\f O$ and \eqref{dva}.
By Theorem \ref{krek},  the map
$\pm T_{d}:\f O_n\rightarrow \f O_n$ is  a minimal holomorphic map between orbifolds, implying that $\chi(\f O_n)\geq 0$. However, it is easy to see that for $n$ big enough this inequality  holds only if $\f O$ is defined either by \eqref{tri}, or by \eqref{chet}, or by  
\be \l{sad} \nu(-1)=\nu(1)=2, \ \ \ \nu(\infty)=n',\ \ \ \ n'\geq 1.\ee 
Finally, checking condition \eqref{uu} at the points of $\pm T_d^{-1}\{-1,1,\infty\}$, we see that  
in the last case $d$ and $n'$ must be coprime. 
 \qed

\bl \l{korova01} Let $A$ be a rational function of degree $d\geq 2$ such that some iterate $A^{\circ l}$ is 
conjugate to $z^{\pm dl}$. Then $A$ is conjugate to $z^{\pm d}$. Similarly, if  
 $A^{\circ l}$ is 
conjugate to $\pm T_{dl}$, then $A$ is conjugate to $\pm T_{d}$.
\el 
\pr 
Assume say that $A^{\circ l}$ is 
conjugate to $z^{\pm dl}$. Then  for any $n$ coprime with $dl$ there exists an orbifold $\f O$ with the signature 
$\{n,n\}$  such that  $z^{\pm dl}\in \f E(\f O),$ implying by Corollary \ref{lat} that $A \in \f E(\f O).$
It follows now from Corollary \ref{l2} that   $A$ is conjugate to $z^{\pm d}$. The case 
where $A^{\circ l}$ is 
conjugate to $\pm T_{dl}$ is considered similarly. \qed



\vskip 0.2cm

\noindent{\it Proof of Theorem \ref{uni}.}
In order to prove the existence of $\f O_0^A$ it is enough to show that there exist at most finitely many orbifolds $\f O$
such that $A:\, \f O\rightarrow \f O$ is a minimal holomorphic map.
Indeed, it follows from Theorem \ref{krek} that in this case we can set 
$$\f O_0^A=\LCM\Big(\f O_1,\f O_2, \dots ,\f O_l\Big),$$ where $\f O_1$, $\f O_2, \dots \f O_l$ is a complete list of such orbifolds.

Assume in contrary that there exists  an infinite sequence of pairwise distinct orbifolds $\f O_1,$ $\f O_2,$ $\dots $   such that $A:\, \f O_i\rightarrow \f O_i$ 
is a minimal holomorphic map for every $i\geq 0.$ Set
$$\f U_s=\LCM\Big(\f O_1,\f O_2, \dots ,\f O_s\Big),\ \ \ s\geq 1.$$ By Theorem \ref{krek},
the maps  $A:\, \f U_s\rightarrow \f U_s$, $s\geq 1,$ are minimal holomorphic maps between orbifolds. Clearly, if the set
$\f U_1,$ $\f U_2,$ $\dots $ 
 is finite, then the set  $\f O_1,$ $\f O_2,$ $\dots $ is also finite. Therefore, the set
$\f U_1,$ $\f U_2,$ $\dots $ 
 is infinite.
Since  $\chi(\f U_s)\geq 0$  by Proposition \ref{p1}, and  $\f U_s\preceq \f U_{s+1}$, this implies that 
for $s$ big enough either $\nu(\f U_s)=\{n,n\},$ or $\nu(\f U_s)=\{2,2,n\},$
where $n\to \infty$ as $s\to \infty.$
However, in this case Corollary \ref{l2} and Corollary \ref{l2+}  imply that the function $A$ is conjugate either to $z^{\pm d}$ or to $\pm T_d,$
in contradiction with the assumption.

By Lemma \ref{korova01}, the orbifolds $\f O_0^{A^{\circ l}}$ either exist for all $l\geq 1$, 
or do not exist  for all $l\geq 1$. Assuming that they exist, 
the proof of the equality 
\be \l{aqaq} \f O_0^{A^{\circ l}}=\f O_0^A\ee
is obtained by a modification of the proof of Corollary \ref{lat}. 
Set $$\f O'=A^*\left(\f O_0^{A^{\circ l}}\right), \ \ \ \ \ \t{\f O}=\LCM(\f O_0^{A^{\circ l}},\f O').$$
 Then $A:\, \f O'\rightarrow \f O_0^{A^{\circ l}}$  and $A^{\circ l}:\, \t{\f O}\rightarrow \t{\f O}$ are  minimal holomorphic maps. 
Since $\f O_0^{A^{\circ l}}\preceq\t{\f O}$, it follows from 
the maximality  of $\f O_0^{A^{\circ l}}$ that 
$\t{\f O}=\f O_0^{A^{\circ l}}.$ This condition is stronger than the condition $\chi(\t{\f O})\geq 0$ used in  Corollary \ref{lat} and combined with 
$\nu(\f O')=\nu(\f O_0^{A^{\circ l}})$ 
 implies that $\f O'=\f O_0^{A^{\circ l}}$. Thus, $A:\, \f O_0^{A^{\circ l}}\rightarrow \f O_0^{A^{\circ l}}$ 
is a minimal holomorphic map, and hence 
$\f O_0^{A^{\circ l}}\preceq\f O_0^{A}$. On the other hand, the first part of Theorem \ref{comp} implies that $\f O_0^{A}\preceq\f O_0^{A^{\circ l}}$. Therefore, \eqref{aqaq} holds.
 \qed

\vskip 0.2cm

Notice that generalized Latt\`es maps are exactly rational functions for which the orbifold $\f O_0^A$ is distinct from the non-ramified sphere, completed
by the functions conjugate to $z^{\pm d}$ or $\pm T_d$ for which the orbifold $\f O_0^A$ is not defined. Furthermore, 
the following statement holds. 

\bl A rational function is a  Latt\`es map if and only if  $\chi(\f O_0^A)= 0.$
\el
\pr The ``if'' part is obvious. On the other hand, if  
 $A:\f O\rightarrow \f O$ is a  covering map, then it follows from 
 $\f O\preceq \f O_0^A$ that $\chi(\f O)\geq \chi(\f O_0^A)$.
Therefore, since \linebreak $\chi(\f O_0^A)\geq 0$ and $\chi(\f O)=0$, the equality  $\chi(\f O_0^A)=0$ holds.
\qed


\br\l{xeretz} \normalfont The functions $z^{\pm n}$ and $\pm T_n$ can be considered as {\it covering} self-maps between orbifolds if to allow the base Riemann surface to be 
non-compact. Namely, it is easy to see that  the map $z^{\pm n}: \f O\rightarrow \f O$ is a covering map for the non-ramified orbifold with the base surface $\f R=\C\setminus\{0,\infty\}$, while 
$\pm T_n: \f O\rightarrow \f O$ is a covering map for the orbifold defined on $\f R=\C\setminus \{\infty\}$ by the condition $\nu(1)=2,$ $\nu(-1)=2$. The corresponding functions $\theta_{\f O}$ are $e^z$ and $\cos z$.  
Notice that the functions $z^{\pm n}$ and $\pm T_n$ along with Latt\`es maps play a key role in the description of commuting 
rational functions obtained by Ritt (see \cite{r}, \cite{e2}, \cite{pnew}). 
\er

In order to check whether or not a given rational function $A$ is a generalized Latt\`es map
one can use the
following lemma. 

\bl \l{toch} Let $A$ be a rational function of degree at least five, and $\f O_1$, $\f O_2$ orbifolds distinct from the non-ramified sphere such that  $A:\f O_1\rightarrow \f  O_2$ is a minimal holomorphic map between orbifolds. Assume that $\chi(\f O_1)\geq 0$. Then    $c(\f O_2)\subseteq c(\f O_2^A)$. \qed
\el 
\pr Suppose that  $z_0\in c(\f O_2)$ is not a critical value of $A.$ Then \eqref{rys} implies that for every point $z\in A^{-1}\{z_0\}$ we have 
$\nu_1(z)=\nu_2(z_0)>1$, implying that $c(\f O_1)$ contains at least five points in contradiction with $\chi(\f O_1)\geq 0.$ 
\qed
\bc \l{sosk}
Let $A$ be a rational function of degree at least five, and $\f O$ an orbifold distinct from the non-ramified sphere such that $A:\f O\rightarrow \f  O$ is a minimal holomorphic map between orbifolds. Then  $c(\f O)\subseteq c(\f O_2^A)$. \qed
\ec

\vskip 0.2cm 
Corollary \ref{sosk} provides a practical method for 
 finding  $\f O_0^A$.
Indeed, it implies  that for a given rational function $A$ of degree at least five, not conjugate to $z^{\pm d}$ or $\pm T_{d}$, any orbifold $\f O$ such that \eqref{uu} holds satisfies $c(\f O)\subseteq c(\f O_2^A)$. Combined with Corollary \ref{l2} and Corollary  \ref{l2+}, this implies that  there exist 
only finitely many possibilities for $\f O$. Finally, for each possible $\f O$ it is enough to check  condition \eqref{uu}  only at points of the finite set \eqref{krysa}.

\end{section}

\begin{section}{Generalized Latt\`es maps for the signatures $\{2,3,3\}$, $\{2,3,4\}$, and $\{2,3,5\}$}

In this section we describe an approach to the description of   minimal holomorphic maps	 $A:\, \f O\rightarrow \f O$
for  $\f O$ with 
$\chi(\f O)>0$ basing on a link between such maps and rational functions $F$ commuting with $\Gamma_{\f O}$. 
We also describe the class of  {\it polynomial}
generalized Latt\`es maps. 
Denote by $Out(\Gamma_{\f O})$ the outer automorphism group of $\Gamma_{\f O}$, and by $d_{\f O}$ the order of 
$Out(\Gamma_{\f O})$.

\bl \l{asz}   Let $\f O$ be an orbifold with $\chi(\f O)>0$,  $A$ a rational function  such that 
$A:\, \f O \rightarrow \f O$ 
is a minimal holomorphic map between orbifolds, and
 $F$ a rational function such  that    diagram \eqref{di} commutes.  
Then there exists  $\sigma\in \Gamma_{\f O}$ such that 
$\sigma\circ F^{\circ d_{\f O}}$ commutes with  $\Gamma_{\f O}$ and the diagram 
\be \l{dr}
\begin{CD}
\C\P^1 @>\sigma \circ F^{\circ d_{\f O}}>> \C\P^1\\
@VV\theta_{\f O}V @VV\theta_{\f O}V\\ 
\C\P^1 @>A^{\circ d_{\f O}} >> \C\P^1 
\end{CD}
\ee
commutes.
\el 
\pr 
Recall that by Proposition \ref{poiu} a rational function  $F$ satisfying \eqref{di} for given $A$ and $\theta_{\f O}$ is defined up to the composition $\sigma\circ F$, where $\sigma \in \Gamma_{\f O}$. 
Furthermore, it is easy to see that for  $\sigma\in \Gamma_{\f O}$ the change 
$F\rightarrow \sigma\circ F$   corresponds to the change $\phi\rightarrow \sigma\circ \phi \circ \sigma^{-1}$. In particular, if the automorphism  $\phi$ is inner, then for an appropriate $\sigma$ 
the automorphism $\sigma\circ \phi \circ \sigma^{-1}$ is identical, or equivalently  
the function  $\sigma\circ F$ commutes with $\Gamma_{\f O}$. Therefore, since \eqref{di} implies 
the equalities $$A^{\circ n}\circ \theta_{\f O}= \theta_{\f O}\circ F^{\circ n}, \ \ \ n\geq 1,$$ 
 $$F^{\circ n}\circ \sigma=\phi^{\circ n}(\sigma)\circ F^{\circ n}, \ \ \ \sigma \in \Gamma_{\f O},
$$ and the automorphism  $\phi^{\circ d_{\f O}}$ is inner, there exists $\sigma\in \Gamma_{\f O}$ 
as required.
\qed

\vskip 0.2cm

Notice that if $\f O$ is given by \eqref{norm}, then
a rational function $F=z^rR(z^n)$
from Theorem \ref{ori1} commutes with $\Gamma_{\f O}=C_n$ if and only if $r=1.$ Thus, since  $d_{\f O}=\phi(n)$, where 
$\phi(n)$ is the Euler totient function, 
the Lemma \ref{asz}  is equivalent in this case to the Euler theorem saying that 
$$r^{\phi(n)}\equiv 1\, \mod n$$
whenever $\GCD(r,n)=1.$ Further,
since $Out(S_4)$ is trivial, Lemma \ref{asz} reduces the description of minimal holomorphic maps  $A:\, \f O\rightarrow \f O$ for orbifolds $\f O$ with  $\nu(\f O)=\{2,3,4 \}$ 
to  the description of rational functions commuting with $S_4.$ 
On the other hand, since $$Out(A_5)=Out(A_4)=\Z/2\Z,$$ it follows from Lemma \ref{asz} that in order to describe 
all minimal holomorphic maps 
$A:\, \f O\rightarrow \f O$ with 
$\nu(\f O)=\{2,3,3\}$ or $\nu(\f O)=\{2,3,5\}$ it is enough to describe the maps corresponding to functions commuting  with $\Gamma_{\f O}$ as well as  
 ``compositional square roots'' of such maps.
The method for describing  rational functions commuting with finite automorphism groups of $\C\P^1$ was given in \cite{dm}. We overview it  below. 

Identify a rational function $f$ with its {\it dual 1-form} as follows. 
Take a representation $f=f_1/f_2,$ where $f_1$ and $f_2$ are polynomials without common roots, construct the homogenization $F_i$ of $f_i$ to the degree $n=\max\{\deg f_1,\deg f_2\}$, and 
set $$\omega=-F_2dx+F_1dy.$$ It is clear that the form $\omega$ is defined up to a multiplication by $\lambda\in \C\setminus\{0\}$, and forms $\omega_1$ and $\omega_2$ represent the same function if and only if $\omega_2=\lambda \omega_1$ for some $\lambda \in \C\setminus\{0\}.$
Under this identification the function 
$\mu^{-1}\circ f\circ \mu,$ 
where $$\mu =\frac{\alpha z+\beta}{\gamma z+\delta}, \ \  \ \alpha,\beta,\gamma,\delta \in \C,$$ 
is identified with the pullback $\mu^{\prime*}\omega,$ where 
$$\mu^{\prime}:\, (x,y)\longrightarrow (\alpha x+\beta y, \gamma x+\delta y).$$
Thus, the problem of describing rational functions commuting with a group $\Gamma$ reduces to the problem of describing forms $\omega$ such that for any $\mu\in \Gamma$ the equality 
$$\mu^{\prime*}\omega=\chi(\mu)\omega,$$
holds for some $\chi(\mu)\in \C.$ On the other hand,  it was shown in \cite{dm}, that a  1-form of degree $n$ satisfies this condition if and only if 
\be \l{form} \omega=U(x,y)\lambda+dV(x,y),\ee where 
$U$ and $V$ are invariant homogeneous polynomials with the same character, $\deg V=n+1,$ $\deg U=n-1$, and 
$$\lambda=-ydx+xdy.$$

It is easy to see that the function $f$ corresponding to form  \eqref{form} is obtained  by setting $z=x/y$ in
\be \l{sho} \frac{xU(x,y)+\frac{\partial V}{\partial y}(x,y)}{yU(x,y)-\frac{\partial V}{\partial x}(x,y)}.\ee
Notice that since $0$ is a form of every  degree, $U$ and $V$ can be equal zero. 
In particular, for any homogeneous polynomial $V$ we obtain a function commuting with $\Gamma$ setting  $z=x/y$ in
\be \l{blik} -\frac{\frac{\partial V}{\partial y}(x,y)}{\frac{\partial V}{\partial x}(x,y)}\,.\ee
On the other hand, if $V=0$, then for any $U$ formula \eqref{sho} leads to the same function  
$f=z$.

Let us illustrate the above considerations by finding explicitly all rational functions of degree $\leq 7$ commuting with the group 
$\Gamma_{\f O}$ for an orbifold $\f O$ with $\nu (\f O)=\{2,3,3\}$, and corresponding 
minimal holomorphic maps $A:\, \f O\rightarrow \f O$.  
According to Klein \cite{klein}, homogenous polynomials for the corresponding group $\Gamma=A_4$ are polynomials in the forms
$$\Phi={x}^{4}+2\,i\sqrt {3}{x}^{2}{y}^{2}+{y}^{4},$$
$$\Psi={x}^{4}-2\,i\sqrt {3}{x}^{2}{y}^{2}+{y}^{4},$$ 
$$t=xy(x^4-y^4).$$
Furthermore, $t$ is absolutely invariant, while $\Phi$ and $\Psi$ are invariant with cha\-racters $\chi_{\Phi}$ and $\chi_{\Psi}$ whose product is the trivial character. 
This implies that all forms \eqref{form} of degree $\leq 6$ are obtained from \eqref{blik} for $V$ equal $\Phi,$ $\Psi,$ or $t.$ Indeed, for non-zero  $U$ and $V$ such a  form 
may satisfy the condition $\deg V=\deg U+2$ only if $U$ is equal to $\Phi$ or $\Psi$, and $V$ is equal to $t$. However, for such $U$ and $V$ the condition concerning characters is not true.
Rational functions commuting with $\Gamma=A_4$ which correspond to forms \eqref{blik} with $V$ equal $\Phi,\Psi,t$  are
$$F_1=-{\frac {i\sqrt {3}{z}^{2}+1}{z \left( i\sqrt {3}+{z}^{2} \right) }},$$
$$F_2=-{\frac {i\sqrt {3}{z}^{2}-1}{z \left( i\sqrt {3}-{z}^{2} \right) }},$$
$$F_3=-{\frac {z \left( {z}^{4}-5 \right) }{5\,{z}^{4}-1}}.$$ 
For the degree seven we obtain a one-parameter series 
setting in \eqref{form}
$$U=ct, \ \ \ \ c\in \C, \ \ \ \ \ V=\Phi\Psi.$$
In order to obtain the corresponding   generalized Latt\`es map in a compact form, it is convenient to rescale this parametrization setting $c=
8i\sqrt{3}a,$  $a\in \C,$ so that  
$$F_4=\frac{1}{z}
\left(
\frac{3\,a{z}^{6}-7\,i{z}^{4}\sqrt {3}-3\,a{z}^{2}-i\sqrt{3}}
{i\sqrt{3}{z}^{6}+3\,a{z}^{4}+7\,i\sqrt {3}{z}^{2}-3\,a}\right).
$$

The generalized Latt\`es maps corresponding to  $F_i$, $1\leq i \leq 4,$  are 
$$L_1={\frac {27 z}{ \left( 4\,z-1 \right) ^{3}}},$$

$$L_2=-{\frac { \left( z-4 \right) ^{3}}{27{z}^{2}}},$$

$$L_3=-{\frac { \left( 5\,z-4 \right) ^{3}}{{z}^{2} \left( 4\,z-5 \right) ^{
3}}},$$
and 
$$L_4=
z \left( {\frac { \left( a-1 \right) ^{4}{z}^{2}-2\, \left( a-1
 \right)  \left( {a}^{3}-3\,{a}^{2}-9\,a-21 \right) z+ \left( a-7
 \right)  \left( a+1 \right) ^{3}}{ \left( a+7
 \right)  \left( a-1 \right) ^{3}{z}^{2}-2\, \left( a+1 \right) 
 \left( {a}^{3}+3\,{a}^{2}-9\,a+21 \right) z+ \left( a+1 \right) ^{4}
 }}\right) ^{3}.
$$
The functions $L_i$, $1\leq i \leq 4,$ and $F_i$, $1\leq i \leq 4,$ are related by the commutative 
diagram
\be 
\begin{CD}
\C\P^1 @>F_i>> \C\P^1 \\
@VV\theta_{\f O} V @VV\theta_{\f O} V\\ 
\C\P^1 @>L_i >> \ \ \C \P^1\, ,
\end{CD}
\ee
where $\f O$ is normalized by the condition
\be \l{orbi} \nu(0)=3, \ \ \ \nu(1)=2, \ \ \ \nu(\infty)=3,\ee and the function 
$$\theta_{\f O}= {\frac { \left( {z}^{4}+2\,i\sqrt {3}{z}^{2}+1 \right) ^{3}}{ \left( {
z}^{4}-2\,i\sqrt {3}{z}^{2}+1 \right) ^{3}}}
$$ is obtained from $\Psi^3/\Phi^3$ by setting $z=x/y.$

Of course, the fact that  $L_i:\, \f O\rightarrow \f O$, $1\leq i \leq 4,$ are indeed minimal holomorphic maps between orbifolds can be checked directly.
 For example, for $L_4$ we must check condition \eqref{uu} at points of the set $L_4^{-1}\{0,1,\infty\}$. 
Clearly, \eqref{uu} holds for any point $z$ such that $L_4(z)=\infty$, since  all  points of $L_4^{-1}\{\infty\}$ distinct from $\infty$ have the multiplicity divisible by $3$ while the multiplicity of $\infty$ is one. Similarly,  \eqref{uu} holds for  points $z$ with $L_4(z)=0.$ Finally, 
formula
 $$\scriptstyle
L_4-1=
(z-1)  \frac{\left((a-1)^6z^3 - \left( 3\,{a}^{3}+3\,{a}^{2}+45\,a+109 \right)\left( a-1 \right) ^
{3}z^2+ \left( 3\,{a}^{3}-3\,{a}^{2}+45\,a-109 \right)  \left( a+1 \right) ^{
3}z-(a+1)^6
\right)^2}
{\left(\left( a+7
 \right)  \left( a-1 \right) ^{3}{z}^{2}-2\, \left( a+1 \right) 
 \left( {a}^{3}+3\,{a}^{2}-9\,a+21 \right) z+ \left( a+1 \right)^{4}\right)^3
 }
$$ implies 
that \eqref{uu} holds  for points $z$ with $L_4(z)=1$.

Notice that the functions $L_1$ and $L_2$ are conjugate by the function $\mu=1/z$.\footnote{We thank to Benjamin Hutz who draw our attention to this fact.} This is explained by the symmetry of the orbifold $\f O$ given by \eqref{orbi} with respect to $\mu$, implying that if $L:\f O\rightarrow \f O$ is a minimal holomorphic map between orbifolds, then $\mu^{-1}\circ L\circ \mu$ is also such a map. Correspondingly, $L_1$ and $L_2$ are conjugate by $\mu$, the function $L_3$ commutes with $\mu,$ and
$$\mu^{-1}\circ L_4(a,z)\circ \mu=L_4(-a,z).$$

In conclusion, we describe the class of  {\it polynomial}
generalized Latt\`es maps. 


\bt \l{polyn} Let $A$ be  a polynomial of degree at least two such that  $A\,: \f O\rightarrow \f O$   is a minimal holo\-mor\-phic  map between 
orbifolds for some  $\f O$ distinct from the non-ramified sphere. Then either 
$A$ is conjugate to $z^rR^n(z)$, where $R\in \C[z]$  and $\GCD(r,n)=1,$ or  
$A$ is conjugate to  $\pm T_m$,  where  $\GCD(m,n)=1.$ 
\et 
\pr Show first that $\chi(\f O)>0.$ Indeed, if $\chi(\f O)=0,$ then  
arguing as in the proof of Theorem \ref{bas}
we can construct commutative diagram \eqref{xui} with $g(R)=1.$
Since $A$ is a polynomial, $A^{-1}\{\infty\}=\infty$, implying that  the set 
$S=\pi^{-1}\{\infty\}$ is completely invariant with respect to $B$.
On the other hand, since $g(R)=1$,  the map $B$ is non-ramified by the Riemann-Hurwitz formula, implying that the set
$B^{-1}(S)$ contains $$\vert S\vert\,\deg B\geq 2\vert S\vert>\vert S\vert$$ points.

Assume now that $\chi(\f O)>0$, and consider diagram \eqref{di} provided by \linebreak Theo\-rem \ref{las}.  It is well known that  the complete $F$-invariance of a finite set implies that it contains at most two points. Therefore, the set $S=\theta_{\f O}^{-1}\{\infty\}$ contains at most two points, 
 and without loss of generality we  
may assume that either $S=\{\infty\}$, or $S=\{0,\infty\}$. 
Since $\theta_{\f O}^{-1}\{\infty\}$ is an orbit 
of $\Gamma_{\f O}$, where $\Gamma_{\f O}$ is one of the five finite rotation groups of the sphere, 
in the first case it follows from $\vert\theta_{\f O}^{-1}\{\infty\}\vert=1$ that $\nu(\f O)=\{n,n\}$, $n\geq 2$. Therefore, since $\theta_{\f O}^{-1}\{\infty\}=\{\infty\},$ without loss of generality  
we may assume that 
$\theta_{\f O}=cz^n$, $c\in \C.$ Moreover, considering instead of the polynomial $A$  
the polynomial $A(cz)/c$, we can assume that $\theta_{\f O}=z^n$. Arguing now as in the proof of Theorem \ref{ori1} and taking into account  that $F$ is a polynomial since $F^{-1}\{\infty\}=\{\infty\}$, 
we obtain that 
 $A=z^rR^n(z)$, where $R$ is a  polynomial and     $\GCD(r,n)=1.$

Similarly, if $S=\{0,\infty\}$, then it follows from $\vert\theta_{\f O}^{-1}\{\infty\}\vert=2$ that
 without loss of generality  
we may assume that \be \l{ser} \theta_{\f O}=\mu\circ \frac{1}{2}\left(z^n+\frac{1}{z^n}\right), \ \ \ n\geq 1,\ee
for some M\"obius transformation $\mu$ such that $\mu(\infty)=\infty.$ 
Indeed, if $\theta_{\f O}^{-1}\{\infty\}$ is a singular orbit 
of $\Gamma_{\f O}$, then  $\nu(\f O)=\{2,2,n\}$, $n\geq 2,$ implying \eqref{ser}. 
On the other hand, if the orbit $\theta_{\f O}^{-1}\{\infty\}$ is non-singular, then  $\nu(\f O)=\{2,2\}.$ 
Therefore, 
$$\theta_{\f O}=az+\frac{b}{z}+c$$ for some $a,b,c\in \C,$ implying that  
composing  $\theta_{\f O}$ with $\sqrt{b/a}\,z$ we still can assume that \eqref{ser} holds. 
Moreover,  since $\mu(\infty)=\infty$, the transformation  $\mu$ is a polynomial, so conjugating $A$ by $\mu$ we can assume that $\mu$ is the identical mapping.

The equality $F^{-1}\{0,\infty\}=\{0,\infty\}$ implies that 
$ F=cz^{\pm m},$ $ c\in \C.$ 
On the other hand,  by Theorem \ref{las}, the homomorphism $\phi$ in \eqref{di+} is an automorphism 
implying that $F$ injectively maps any fiber of $\theta_{\f O}$ onto another fiber. Therefore, the singular 
fiber $ \theta_{\f O}^{-1}\{1\}$ consisting of $n$th roots of 1 is mapped either to itself or to the other singular fiber 
$\theta_{\f O}^{-1}\{-1\}$ consisting of $n$th roots of -1. 
Since this implies that $c^2$ is an $n$th root of unity, 
it follows now from 
$$ A\circ  \frac{1}{2}\left(z^n+\frac{1}{z^n}\right)= \frac{1}{2}\left(z^n+\frac{1}{z^n}\right)\circ cz^d=\pm T_d\circ \frac{1}{2}\left(z^n+\frac{1}{z^n}\right)$$
that $A$ is conjugate to $\pm T_d.$ \qed

\end{section}

\bibliographystyle{amsplain}

\end{document}